\documentclass[11pt,leqno,%dvipdfmx
]{amsart}
\usepackage{amsmath,epsfig,graphicx,color}
\usepackage{comment}
\usepackage{ascmac}
\definecolor{darkblue}{rgb}{.2, 0.2,.8}
\definecolor{carageen}{rgb}{0,0.5,0.3}
\definecolor{darkred}{rgb}{1, 0,0}
\newcommand{\bfxi}{{\boldsymbol \xi}}

\renewcommand{\a}{\alpha}
\newcommand{\red}{\color{darkred}}
\newcommand{\blue}{\color{darkblue}}
\newcommand{\green}{\color{carageen}}

\newcommand{\clt}{central limit theorem}

\newcommand{\ex}{{\rm e}\,}

\newcommand{\asy}{asymptotic}

\newcommand{\ts}{time series}

\newtheorem{lemma}{Lemma}[section]

\newtheorem{theorem}[lemma]{Theorem}

\renewcommand{\P}{{\mathbb P}}
\newtheorem{proposition}[lemma]{Proposition}
\newtheorem{definition}[lemma]{Definition}
\newtheorem{corollary}[lemma]{Corollary}
\newtheorem{example}[lemma]{Example}
\newtheorem{exercise}[lemma]{Exercise}
\newtheorem{remark}[lemma]{Remark}

\newtheorem{tab}[lemma]{Table}

\newcommand{\bfQ}{{\bf Q}}
\newcommand{\bfM}{{\bf M}}
\newcommand{\bfu}{{\bf u}}

\newcommand{\bfR}{{\bf R}}

\newcommand{\bth}{\begin{theorem}}
\newcommand{\ethe}{\end{theorem}}

\newcommand{\bre}{\begin{remark}\em }
\newcommand{\ere}{\end{remark}}

\newcommand{\ble}{\begin{lemma}}
\newcommand{\ele}{\end{lemma}}
\newcommand{\sre}{stochastic recurrence equation}
\newcommand{\pp}{point process}
\newcommand{\bde}{\begin{definition}}
\newcommand{\ede}{\end{definition}}
\newcommand{\bco}{\begin{corollary}}
\newcommand{\eco}{\end{corollary}}

\newcommand{\bpr}{\begin{proposition}}
\newcommand{\epr}{\end{proposition}}

\newcommand{\bexer}{\begin{exercise}}
\newcommand{\eexer}{\end{exercise}}

\newcommand{\bexam}{\begin{example}}
\newcommand{\eexam}{\end{example}}

\newcommand{\efi}{\end{fig}}

\newcommand{\btab}{\begin{tab}}
\newcommand{\etab}{\end{tab}}

\newcommand{\lhs}{left-hand side}

\newcommand{\rv}{random variable}

\newcommand{\sign}{{\rm sign}}

\newcommand{\cov}{{\rm cov}}

\newcommand{\bfTh}{\mbox{\boldmath$\Theta$}}

\newcommand{\rhs}{right-hand side}
\newcommand{\df}{distribution function}

\newcommand{\dint}{\displaystyle\int}

\newcommand{\beao}{\begin{eqnarray*}}
\newcommand{\eeao}{\end{eqnarray*}\noindent}

\newcommand{\beam}{\begin{eqnarray}}
\newcommand{\eeam}{\end{eqnarray}\noindent}

\newcommand{\beqq}{\begin{equation}}
\newcommand{\eeqq}{\end{equation}\noindent}

\newcommand{\bce}{\begin{center}}
\newcommand{\ece}{\end{center}}

\newcommand{\barr}{\begin{array}}
\newcommand{\earr}{\end{array}}

\newcommand{\stp}{\stackrel{\P}{\rightarrow}}
\newcommand{\std}{\stackrel{d}{\rightarrow}}

\newcommand{\stv}{\stackrel{v}{\rightarrow}}
\newcommand{\stw}{\stackrel{w}{\rightarrow}}

\newcommand{\eqd}{\stackrel{d}{=}}

\newcommand{\vague}{\stackrel{\lower0.2ex\hbox{$\scriptscriptstyle
                    \it{v} $}}{\rightarrow}}
\newcommand{\weak}{\stackrel{\lower0.2ex\hbox{$\scriptscriptstyle
                    \it{w} $}}{\rightarrow}}
\newcommand{\what}{\stackrel{\lower0.2ex\hbox{$\scriptscriptstyle
                    \it{\hat{w}} $}}{\rightarrow}}

\newcommand{\bdis}{\begin{displaymath}}
\newcommand{\edis}{\end{displaymath}\noindent}

\newcommand{\R}{\mathbb{R}}

\newcommand{\nto}{n\to\infty}
\newcommand{\kto}{k\to\infty}
\newcommand{\xto}{x\to\infty}

\newcommand{\wt}{\widetilde}
\newcommand{\wh}{\widehat}
\newcommand{\vep}{\varepsilon}

\newcommand{\la}{\lambda}

\newcommand{\regvary}{regularly varying}

\newcommand{\regvar}{regular variation}

\newcommand{\bbr}{{\mathbb R}}

\newcommand{\bbz}{{\mathbb Z}}
\newcommand{\Z}{{\mathbb Z}}

\newcommand{\con}{convergence}

\newcommand{\st}{such that}
\newcommand{\fif}{if and only if}
\newcommand{\wrt}{with respect to}
\newcommand{\chf}{characteristic function}
\newcommand{\fct}{function}

\newcommand{\ds}{distribution}

\newcommand{\cmt}{continuous mapping theorem}
\newcommand{\seq}{sequence}

\newcommand{\pro}{probabilit}

\newcommand{\ms}{measure}

\newcommand{\bfx}{{\bf x}}
\newcommand{\bfX}{{\bf X}}

\newcommand{\bfY}{{\bf Y}}
\newcommand{\bfy}{{\bf y}}

\newcommand{\bfZ}{{\bf Z}}

\newcommand{\bfS}{{\bf S}}

\newcommand{\E }{{\mathbb E}}
\renewcommand{\P }{{\mathbb P}}

\newcommand{\1}{{\mathbf 1}}

\allowdisplaybreaks

\begin{document}
\today
\bibliographystyle{plain}
\title[Self-normalized partial sums of heavy-tailed time series]{Self-normalized partial sums of heavy-tailed time series}

\thanks{Muneya Matsui's research is partly supported by the JSPS Grant-in-Aid for Scientific Research C
(19K11868). Thomas Mikosch's research is partly supported by the grant No
9040-00086B
of the Danish Free Research Council (DFF)} 

\author[M. Matsui]{Muneya Matsui}
\address{Department of Business Administration, Nanzan University, 18
Yamazato-cho, Showa-ku, Nagoya 466-8673, Japan.}
\email{mmuneya@gmail.com}

\author[T. Mikosch]{Thomas Mikosch}
\address{Department  of Mathematics,
University of Copenhagen,
Universitetsparken 5,
DK-2100 Copenhagen,
Denmark}
\email{mikosch@math.ku.dk}
\author[O. Wintenberger]{Olivier Wintenberger}
\address{Sorbonne Universit\'e,
4 Place Jussieu, 
75005 Paris,
France}
\email{olivier.wintenberger@sorbonne-universite.fr} 

\begin{abstract}
We study the joint limit behavior of sums, maxima and $\ell^p$-type moduli for samples taken from an $\bbr^d$-valued  \regvary\ stationary \seq\ with infinite variance. As  a con\seq , we can determine the \ds al limits for ratios of sums and maxima, studentized sums, and other self-normalized quantities in terms of hybrid \chf s  and Laplace transforms. These transforms
enable one to calculate moments of the limits and to characterize the differences between the iid and stationary cases in terms of indices which describe effects of extremal clustering
on \fct als acting on the dependent \seq .
\end{abstract}
\maketitle

\section{Introduction}\setcounter{equation}{0}
Consider an $\bbr^d$-valued (strictly) stationary \seq\ $(\bfX_t)$ 
with generic element $\bfX$. We call this \ts\ heavy-tailed if
it does not have finite moments of a certain order. To be more 
specific, we are interested in {\em \regvary\ \ts }.
This notion was introduced by Davis and Hsing \cite{davis:hsing:1995}
via the vague \con\ of certain tail \ms s. We will stick to an equivalent
formulation via weak \con\ conditional on a large value of $|\bfX_0|$; see
Basrak and Segers \cite{basrak:segers:2009}: the stationary \seq\ $(\bfX_t)$
is {\em \regvary\ with (tail) index $\a>0$} if there exist a 
Pareto$(\a)$-distributed \rv\ $Y$, i.e., $\P(Y>y)=y^{-\a}$, $y>1$, 
and an $\bbr^d$-valued \seq\ $(\bfTh_t)$ which is independent of $Y$ \st\
for every $h\ge 0$,
\beao
\P\big(x^{-1}(\bfX_{-h},\ldots,\bfX_h)\in\cdot \,\big|\,|\bfX_0|>x\big)
\std \P\big(Y\,(\bfTh_{-h},\ldots,\bfTh_h)\in \cdot\big)\,,\qquad \xto\,.
\eeao
We call $(\bfTh_t)$ the {\em spectral tail process of $(\bfX_t)$}.
Under suitable {\em mixing} and {\em anti-clustering conditions}
one has \con\ of the 
normalized maxima 
\beao
a_n^{-1}M^{|\bfX|}_n=a_n^{-1}\max_{i=1,\ldots,n}|\bfX_i|\std \eta_\a\,,\qquad \nto\,, 
\eeao
where $\eta_\a$ has a scaled Fr\'echet \ds , i.e., 
$\P(\eta_\a\le x)=\Phi_\a^{\theta_{|\bfX|}}(x)=\ex^{-\theta_{|\bfX|}\,x^{-\a}}$, $x>0$,
and the scaling factor $\theta_{|\bfX|}\in (0,1]$ is the  {\em 
extremal index} of the \seq\ $(|\bfX_t|)$. The \seq\ $(a_n)$ is chosen
\st\ $n\,\P(|\bfX|>a_n)\to 1$ as $\nto$.
Moreover, for $\a\in (0,2)$ 
the (centered) partial sums of the \ts\
converge to an $\a$-stable limit $\bfxi_\a$:
\beao
a_n^{-1}\bfS_n=a_n^{-1} \big(\bfX_1+\cdots +\bfX_n\big)\std \bfxi_\a\,,\qquad \nto.
\eeao
We refer to Section~\ref{sec:prelim} for precise definitions and conditions. 
\par
The normalizing \seq\ $(a_n)$ 
is in general unknown and has to be estimated from data.
A way of replacing $(a_n)$ by a known quantity is 
the principle of {\em self-normalization}, including the well-known 
{\em studentization}; see de la Pe\~{n}a et al. \cite{pena:lai:shao:2009}
for a general treatment of self-normalizing techniques.
A natural choice of self-normalization for the partial 
sums $(\bfS_n)$ is given by the \seq\ of the partial 
maxima $(M_n^{|\bfX|})$. Indeed, for $\a\in (0,2)$ both \seq s require
the same normalization $(a_n)$. Thus, if one can show the joint
\con\ $a_n^{-1} (M_n^{|\bfX|},\bfS_n)\std (\eta_\a,\bfxi_\a)$ it follows that
\beao
\dfrac{\bfS_n}{M_n^{|\bfX|}}\std \dfrac{\bfxi_\a}{\eta_\a}\,, \qquad \nto\,,
\eeao
where $\eta_\a$ and $\bfxi_\a$ are dependent even if $(\bfX_t)$ is an iid \seq .
Moreover, we will extend this result to more general self-normalizing \seq s.
In particular, we will show the \con\ of {\em studentized} sums
\beao
\dfrac{\bfS_n}{\gamma_{n,p}}\std \dfrac{\bfxi_\a}{\zeta_{\a,p}}\,,\qquad \nto\,,
\eeao
where 
\beam\label{eq:moduli}
\gamma_{n,p}=\big(\sum_{t=1}^n|\bfX_t|^p\big)^{1/p}\,,\qquad p>\alpha\,,
\eeam
is an 
$\ell^p$-norm-type modulus of the sample
$\bfX_1,\ldots,\bfX_n$, and $\zeta_{\a,p}^p$ is a positive $\a/p$-stable \rv .
In contrast to the $\a$-stable limit $\bfxi_\a$ of $(a_n^{-1}\bfS_n)$ the limit ratios 
$\bfxi_\a/\eta_\a$ and $\bfxi_\a/\zeta_{\a,p}$ have moments of order $\a$.
\par
For $\a>2$, $\bfX$ has finite second moment. Then, under general
weak dependence conditions the (centered) partial sums  $(\bfS_n/\sqrt n)$ satisfy a \clt\
with Gaussian limit while $(a_n^{-1} M_n^{|\bfX|})$ still has a Fr\'echet limit $\eta_\a$, and $(a_n/\sqrt{n})( \bfS_n/M_n^{|\bfX|})$ converges to a ratio
of a Gaussian vector and $\eta_\a$ which, typically, are independent; see  
Anderson and Turkman \cite{anderson:turkman:1991a,anderson:turkman:1991b,anderson:turkman:1995}, Hsing \cite{hsing:1995}. In this case, normalization of 
$\bfS_n$ by $M_n^{|\bfX|}$ is less meaningful because one needs to know $(a_n)$.
\par
The paper is organized as follows. In Section~\ref{sec:prelim}
we introduce the basic conditions used throughout the paper: mixing and anti-clustering.
We discuss these and also properties of a \regvary\ stationary \seq. In Section~\ref{sec:maxsum} we continue with the main
results of this paper, dealing with the joint \con\ of sums and maxima (Theorem~\ref{th:joint}),
and sums, maxima and $\ell^p$-type moduli (Theorem~\ref{pr:prophybridch}). In Section~\ref{sec:ratiolmit}
we apply these results to ratios of sums and maxima, studentized sums, Greenwood statistics, and 
ratios of norms. In Section~\ref{sec:examples} we indicate how the conditions of the results
in the previous sections can be verified. We focus on the class of \regvary\ iterated random Lipschitz \fct s,  including autoregressive processes and solutions to affine \sre s.
\section{Preliminaries}\label{sec:prelim}\setcounter{equation}{0}
\subsection{Mixing condition}
Mixing conditions are standard for proving \asy\ theory for 
sums and maxima of a stationary \seq . Since we are interested in the joint \con\ of these objects
we introduce a mixing condition which is tailored for this situation.
Following Chow and Teugels \cite{chow:teugels:1978}, we will express the mixing 
condition through a combination of a \chf\ (for sums)  and a \df\ (for maxima)
(so-called {\em hybrid \chf )}:  
\beao
\Psi_n(\bfu,x)&:=&\E\big[\exp\big(i\,a_n^{-1} \bfu^\top\bfS_n\big)\,\1\big(a_n^{-1}M_n^{|\bfX|}\le x\big)\big]\,,\qquad \bfu\in\bbr^d\,,\qquad x>0\,. 
\eeao
The hybrid \chf\ of a pair $(\bfY,U)$
of a random vector $\bfY$ and a  non-negative \rv\ $U$ determines the \ds\ of $(\bfY,U)$. Moreover,
$(\bfY_n,U_n)\std (\bfY,U)$ \fif\ the hybrid \chf s of $(\bfY_n,U_n)$ converge
point-wise to the corresponding hybrid \chf\ of $(\bfY,U)$.
\par
We will use the following {\bf mixing condition}:
for some integer \seq s $r_n\to \infty$, $k_n=[n/r_n]\to\infty$,
\beam\label{eq:wdepms}
 \Psi_n(\bfu,x)&=&
\Big(\E\big[\exp\big(i\,a_n^{-1}\bfu^\top\bfS_{r_n}\big)\,\1\big(a_n^{-1}M_{r_n}^{|\bfX|}\le x\big)\big]\Big)^{k_n}+o(1) \,,\qquad \nto\,,\nonumber\\&&\qquad x>0\,,\qquad \bfu\in\bbr^d\,.
\eeam
Condition \eqref{eq:wdepms} ensures the \asy\ independence of $k_n$ maxima
and sums over disjoint blocks of length $r_n$. It
follows from strong mixing properties of $(\bfX_n)$ or by coupling arguments; see for example Rio \cite{rio:2017} as a general reference  and Section~\ref{sec:examples} below.

\subsection{The anti-clustering condition}\label{subsec:anticlclt}
Anti-clustering conditions are also standard for proving limit theory for
sums and maxima of stationary \seq s; see for example Basrak and Segers 
\cite{basrak:segers:2009}. They ensure that clusters of extreme events
cannot last forever.
The following condition was introduced
in Bartkiewicz et al. \cite{bartkiewicz:jakubowski:mikosch:wintenberger:2011};
it is particularly tailored for partial sums.
\par
An $\bbr^d$-valued stationary \regvary\ \seq\ $(\bfX_t)$
satisfies the  {\bf anti-clustering condition} if for some $r_n\to\infty$ \st\ $k_n=[n/r_n]\to\infty$,
\beam\label{cond:acac}
\lim_{k\to\infty}\limsup_{\nto}  n \,\sum_{j=k}^{r_n}\E\big[(|a_n^{-1}\bfX_{j} |\wedge  x)\,(|  a_n^{-1}\bfX_{0}|\wedge
x)\big]=0\,,\qquad x=1\,.
\eeam
Here and throughout the paper
$(a_n)$ satisfies $n\,\P(|\bfX|>a_n)\to 1$ as $\nto$.
Condition \eqref{cond:acac} can be checked by coupling arguments; see for example Kulik et al. \cite{kulik:soulier:wintenberger:2019} and Section~\ref{sec:examples} below.
For $a\ge 0$ we have
$a\wedge x\le x\,(a\wedge1)\,\1(x\ge 1)+
 (a\wedge1)\, \1(x<1)$.
Therefore \eqref{cond:acac} holds for every $x>0$ if it does for $x=1$. 
\par
Condition \eqref{cond:acac} is symmetric in time: since $(\bfX_t)$ 
is stationary it is equivalent to 
\beao
\lim_{k\to\infty}\limsup_{\nto}  n \,\sum_{j=k}^{r_n}\E\big[(|a_n^{-1}\bfX_{-j} |\wedge  1)\,(|  a_n^{-1}\bfX_{0}|\wedge
1)\big]=0\,.
\eeao
Since 
$a\wedge b = a \,\1(a\le b)+ b\,\1(a>b)$, $a,b>0$,
 \eqref{cond:acac} implies
\beao
\lim_{k\to \infty}\limsup_{\nto}  n \,\sum_{j=k}^{r_n}\P\big(|\bfX_{j} |>  x\,a_n \,, \,  |\bfX_{0}|>x\,a_n \big) =0\,,\qquad x>0\,,
\eeao
hence the anti-clustering 
condition
\beao
\lim_{k\to \infty}\limsup_{\nto}   \,\P\big(\max_{j=k,\ldots,r_n}|\bfX_{j} |>  x\,a_n \,\big|\, |\bfX_{0}|>x\,a_n\big) =0\,,\qquad x>0\,,
\eeao
which is standard for proving limit theory for the extremes of a stationary
\seq ; see for example Davis and Hsing \cite{davis:hsing:1995} and Basrak and Segers \cite{basrak:segers:2009}. 
\par
If $\a\in (1,2)$ we have $\E[a_n^{-1}|\bfX|\wedge 1]\le c/a_n$ (here and in what follows, $c$ denotes any positive constant whose value is not of interest)
while for $\a\in (0,1)$ by Karamata's theorem for \regvary\ \fct s (see Bingham
et al. \cite{bingham:goldie:teugels:1987})  
$\E[a_n^{-1}|\bfX|\wedge 1]\le c/n$.
Therefore
\beao
n\,r_n (\E[a_n^{-1}|\bfX|\wedge 1])^2&=&\left\{\barr{ll}O(r_n/n)=o(1)\,, &\mbox{if $\a\in (0,1)$}\,,\\[2mm]O(r_n\,n/a_n^2)=o(1)\,,
&\mbox{if $\a\in (1,2)$ and also $r_n=o(a_n^{2}/n)$.}
\earr\right.
\eeao
Under these conditions on $(r_n)$, \eqref{cond:acac} holds for an 
iid \seq\ $(\bfX_t)$. Moreover, \eqref{cond:acac} turns into
\beam
\label{cond:centeracac} 
\lim_{k\to\infty}\limsup_{\nto}n \,\sum_{j=k}^{r_n}\cov\big(|a_n^{-1}\bfX_{j} |\wedge  1\,,|  a_n^{-1}\bfX_{0}|\wedge
1\big)=0 \,.
\eeam
\subsection{Properties of \regvary\ stationary \seq s}
Write  $\bfTh=(\bfTh_t)_{t\in\bbz}$ for the spectral tail process of $(\bfX_t)$.
Under \eqref{eq:wdepms} and \eqref{cond:acac} we have the property that
$\|\bfTh\|_\a^\a=\sum_{t\in\bbz}|\bfTh_t|^\a<\infty$ a.s.; see 
Janssen \cite{janssen:2019} and Buritic\'a et al. \cite{buritica:meyer:mikosch:wintenberger:2021}. Then one can define
the {\em spectral cluster process} 
$\bfQ=\bfTh/\|\bfTh\|_\a$.  We will also make use of a change of \ms\ of $\bfQ$ in
$\ell^\a(\bbr^d)$ given  by
\beam\label{eq:qtilde}
\P(\wt \bfQ \in \cdot)&=&%\E\Big[\dfrac{|\bfQ_{T ^\ast}|^\a}{\E[|\bfQ_{T ^\ast}|^\a]}\1\Big( \dfrac{\bfQ}{|\bfQ_{T ^\ast}|} \in \cdot\Big)\Big]\nonumber\\
%&=& \dfrac{1}{\theta_{|\bfX|}}\E\Big[|\bfQ_{T ^\ast}|^\a\,\1\Big( \dfrac{\bfQ}{|\bfQ_{T ^\ast}|} \in \cdot\Big)\Big]\\
%&=&  
\P\Big( \dfrac{\bfQ}{\max_{t\in \bbz}|\bfQ_t|} \in \cdot\Big| \max_{t\in\bbz}|Y \bfQ_t|>1\Big)\,,
%&=&\red \E\Big[\dfrac{|\bfTh_{T ^\ast}|^\a}{\E[|\bfTh_{T ^\ast}|^\a]}\1\Big( \dfrac{\bfTh}{|\bfTh_{T ^\ast}|} \in \cdot\Big)\Big]\nonumber\,.
\eeam
and the extremal index of $(|\bfX_t|)$
can be expressed as  
\beam\label {eq:august28b}
\theta_{|\bfX|}=\P(Y\,\max_{t\in\bbz} |\bfQ_t|>1)= \E\big[\max_{t\in\bbz} |\bfQ_t|^\a\big]\,.
\eeam
An alternative way of defining \regvar\ of the stationary process
$(\bfX_t)$ is via the vague \con\ relations on $\bbr^{d(h+1)}_{\bf0}=\bbr^{d(h+1)}\backslash \{\bf0\}$, for $h\ge 0$,
\beam\label{eq:vagconv}
n\,\P\big(a_n^{-1}(\bfX_0,\ldots,\bfX_h)\in\cdot \big)\stv \mu_h(\cdot)\,,\qquad \nto\,,
\eeam
where the {\em tail \ms s} $(\mu_h)$ are non-null and have the homogeneity
property $\mu_h(t\,\cdot)=t^{-\a}\,\mu_h(\cdot)$, $t>0$. We observe that, due to 
stationarity of $(\bfX_t)$, 
these \ms s have the consistency property $\mu_{h+1}(\bbr^d\times \cdot)
=\mu_h(\cdot)= \mu_{h+1}(\cdot\times  \bbr^d)$; see Kulik and Soulier \cite{kulik:soulier:2020}.
\par
An important tool for dealing with \regvary\ stationary \seq s is the {\em time-change
formula} (see Basrak and Segers \cite{basrak:segers:2009}), for $h\ge 0$, $
t \in\bbz$,
\beam\label{eq:tcp}
\P((\bfTh_{-h},\ldots,\bfTh_{h})\in \cdot \mid \bfTh_{-t}\ne {\bf0})&=&
\E \Big[\dfrac{|\bfTh_{t}|^\a}{\E\big[|\bfTh_{t}|^\a\big]} 
\1\Big(\dfrac{(\bfTh_{t-h},\ldots,\bfTh_{t+h})}{|\bfTh_{t}|}\in \cdot\Big)\Big]\,.
\eeam

\section{Joint convergence of sums, maxima and norms}\label{sec:maxsum}\setcounter{equation}{0}
\subsection{Joint convergence of sums and maxima}
In this section we prove our main result about the joint \con\ of normalized
maxima and sums based on a \regvary\ stationary \seq . We restrict ourselves to
indices $\a\in (0,2)\backslash\{1\}$. The main reason for this is that
centering of sums for $\a=1$ is complicated and stretches the proofs due to much
technical detail. Here and in the subsequent sections we use transform 
techniques (hybrid \chf s, Laplace transforms). This is in contrast to
some of the more recent developments in limit theory for \regvary\ stationary \seq s;
various authors prefer to use \pp\ techniques (first proving the weak \con\
of the processes of the points $(a_n^{-1}\bfX_t)$, then applying a.s. continuous
mappings to these points); see for example Davis and Hsing \cite{davis:hsing:1995},
Resnick \cite{resnick:2007}, Kulik and Soulier \cite{kulik:soulier:2020},
Krizmani\'c \cite{krizmanic:2019}. They obtain the series representation ($0<\a<1$) of Lepage et al. \cite{lepage:woodroofe:zinn:1981}
\beam\label{eq:april10a}
a_n^{-1}(M_n^{|\bfX|},\bfS_n)\std \Big(\sup_{i\ge 1}\Gamma_i^{-1/\a}\sup_{j\in\Z}|\bfQ_{ij}|,\sum_{i\ge 1}\Gamma_i^{-1/\a}\sum_{j\in\Z}\bfQ_{ij}\Big)\,,\qquad \nto\,.
\eeam
 When dealing with partial sums for $\a\in (1,2)$, a difficulty 
is the {\em vanishing-small-values condition}, i.e., the sums of the truncated quantities
$(\bfX_t/a_n)$ have to be negligible. We avoid these conditions. A further advantage 
of transform techniques is that it is easier to identify the (joint) limit \ds\ of
sums, maxima, self-normalized sums, and their \ds al characteristics such as moments.
\bth\label{th:joint} 
Consider an $\bbr^d$-valued \regvary\ stationary process $(\bfX_t)$ 
with index $\alpha\in  (0,2)\backslash \{1\}$. If $\a\in (1,2)$ we also assume 
that $\E[\bfX]=0$.
Choose the normalizing constants
$(a_n)$ \st\
$n\,\P(|\bfX|> a_n)\to 1$ as $\nto$.  
Assume the mixing condition \eqref{eq:wdepms} and 
the anti-clustering condition \eqref{cond:acac}
for the same integer \seq s
$r_n\to\infty$, $k_n=[n/r_n]\to\infty$ as $\nto$.
Then 
\beao
a_n^{-1}(M_n^{|\bfX|},\bfS_n)\std (\eta_\a,\bfxi_\a)\,,\qquad \nto\,,
\eeao
where $\eta_\a$ is Fr\'echet-distributed with a positive 
extremal index $\theta_{|\bfX|}$ and $\bfxi_\a$ is $\a$-stable 
with \chf\
\beam\label{eq:agree}
\varphi_{\bfxi_\alpha}(\bfu) = \exp\big( -c_\a\,\sigma^\a(\bfu)\big(1- i\, \beta(\bfu) \tan(\a\,\pi/2) \big) \big) \,,\qquad \bfu\in\bbr^d\,,
\eeam
constant  $c_\alpha = \Gamma(2-\a)\,\cos(\a\,\pi/2)\,/(1-\alpha)$, 
and parameter \fct s 
\beam\label{eq:agree1}
\beta(\bfu)&=&\dfrac{\E\big[(\bfu^\top\sum_{t\in \bbz}\bfQ_t )_+^\a- (\bfu^\top\sum_{t\in \bbz} \bfQ_t )_-^\a \big]}{\E\big[|\bfu^\top\sum_{t\in \bbz}  \bfQ_t |^\a\big]}\,,\qquad
\sigma^\a(\bfu)=\E\Big[\Big|\bfu^\top\sum_{t\in \bbz} \bfQ_t \Big|^\a\Big]\,.%\nonumber 
\eeam
Moreover, 
the joint \ds\ of $(\eta_\a,\bfxi_\a)$ is characterized by the hybrid \chf
\beam\label{eq:maxsum}\lefteqn{
\E\big[\ex^{i \,\bfu^\top\bfxi_\a}\1(\eta_\a\le x)\big]}\nonumber\\
&=&\varphi_{\bfxi_\a}(\bfu)\,
\exp\Big(-\dint_0^\infty\E\Big[\ex^{i\,y\,\bfu^\top \sum_{t=-\infty}^\infty \bfQ_t}\, \1\Big(y\,\max_{t\in \bbz}|\bfQ_t|>x \Big)\Big]\,d(-y^{-\a})\Big)\\ %\nonumber\\
&=&\varphi_{\bfxi_\a}(\bfu)\,\Phi_\a^{\theta_{|\bfX|}}(x)\, \exp\Big(-\theta_{|\bfX|}\dint_x^\infty\E\big[\ex^{i\,y\,\bfu^\top\sum_{t=-\infty}^\infty \wt \bfQ_t}-1\big]\,d(-y^{-\a})\Big)\,,\; \bfu\in\bbr^d\,, x>0\,,\nonumber
\eeam
where $\bfQ=\bfTh/\|\bfTh\|_\a$ and $\wt \bfQ$ is defined in \eqref{eq:qtilde} and the extremal index $\theta_{|\bfX|}$ is given in \eqref{eq:august28b}.
\ethe
In the iid univariate case and for general maximum domains of attraction, 
Chow and Teugels \cite{chow:teugels:1978} showed a 
result of the type \eqref{eq:maxsum}.
\par 
Setting   $x=\infty$ in \eqref{eq:maxsum}, we conclude that 
$a_n^{-1}\bfS_n\std \bfxi_\a$. Similarly, setting $\bfu=\bf0$, we obtain
for $x>0$,
\beao
\P(a_n^{-1}M_n^{|\bfX|}\le x)&\to& \P\big(\eta_\a\le x\big)\\ 
&=&%\varphi_{\bfxi_\a}(\bfu)\,
\exp\Big(-\dint_0^\infty \P\Big(y\,\max_{t\in \bbz}|\bfQ_t|>x \Big)\,d(-y^{-\a})\Big)
\\
&=& \exp\Big(-x^{-\a}\,\E\Big[\max_{t\in \bbz}|\bfQ_t|^\a\Big]\Big)\\
&=&\Phi_\a^{\theta_{|\bfX|}}(x)\,,\qquad \nto\,.
\eeao
In the last step we used the identity \eqref{eq:august28b}. Therefore the factor
\beao
\exp\Big(-\theta_{|\bfX|}\dint_x^\infty\E\big[\ex^{i\,y\,\bfu^\top\sum_{t=-\infty}^\infty \wt \bfQ_t}-1\big]\,d(-y^{-\a})\Big)
\eeao
in the limit hybrid \chf\ \eqref{eq:maxsum} describes the dependence
between $\eta_\a$ and~$\bfxi_\a$. 
\par
In the iid and \asy ally independent cases, i.e.,  when  $\bfTh_t=\bfQ_t=\bf0$, $t\ne 0$, $\bfQ_0=\bfTh_0$, \eqref{eq:maxsum} turns into
\beao
\E\big[\ex^{i\bfu\,\bfxi_\a}\,\1(\eta_\a\le x)\big]&=&\varphi_{\bfxi_\a}(\bfu)\,
\exp\Big(-\dint_x^\infty\E\big[\ex^{i\,y\,\bfu^\top \bfTh_0}\big]\,d(-y^{-\a})\Big)\,,\nonumber\\
&&\bfu\in\bbr^d\,,\qquad x>0\,.
\eeao
The \rhs\ does not factorize into $\varphi_{\bfxi_\a}(\bfu)\P(\eta_\a\le x)$ for all $\bfu\in\bbr^d$ and $x>0$. Therefore $\eta_\a$
and $\bfxi_\a$ are dependent. The hybrid \chf\ $\Psi$ of $(\eta_\a,\bfxi_\a)$ factorizes, i.e., for all $\bfu\in\bbr^d$, $x>0$,
\beao
\Psi(\bfu,x)=
\E\big[\ex^{i\bfu\,\bfxi_\a}\,\1(\eta_\a\le x)\big]=\varphi_{\bfxi_\a}(\bfu)\,
\P(\eta_\a\le x)\,,
\eeao  
\fif\ $\eta_\a$, $\bfxi_\a$ are independent since the quantity 
$\Psi$ determines the  \ds\ of $(\eta_\a,\bfxi_\a)$. In 
view of \eqref{eq:maxsum}, 
$\eta_\a$ and $\bfxi_\a$ are independent \fif\ 
\beao
\dint_x^\infty \E\big[ \ex^{i\,y\,\bfu^\top\sum_{t=-\infty}^\infty \wt \bfQ_t}-1 \big]\,d(-y^{-\a})=0\,,\qquad \bfu\in\bbr^d\,,\qquad x>0\,.
\eeao
If $Y$ is Pareto$(\a)$ and independent of $(\wt \bfQ_t)$ the latter condition
for $x=1$ implies that the \chf\ of $Y \sum_{t=-\infty}^\infty \wt \bfQ_t$ is $1$ which in turn implies that
$\sum_{t\in\bbz}\bfQ_t=\sum_{t\in\bbz}\bfTh_t=\bf0$ a.s. We conclude from the form of
$\varphi_{\bfxi_a}$ that $\bfxi_\a=\bf0$ a.s. Thus, $\eta_\a$ and $\bfxi_\a$ 
are independent \fif\ $\bfxi_\a$ degenerates.

\begin{proof}
It suffices 
to show that
\beam\label{eq:augu10}
\Psi_n(\bfu,x)\to \Psi(\bfu,x)\,,\qquad \nto\,,\qquad  \bfu\in\bbr^d\,,\; x> 0\,.
\eeam
The mixing condition \eqref{eq:wdepms} and a Taylor expansion yield for 
fixed $\bfu\in\bbr^d$, $x>0$, 
\beam\label{eq:xmas8a}
\log\Psi_n(\bfu,x)
%\E\Big[\ex^{i\,u\sum_{t=1}^n X_t/a_n \Big)\1\Big(\max_{1\le t\le n}X_t/a_n>x %\Big)\Big]}\\
&=& k_n\,\big(\E\big[\ex^{i\,a_n^{-1}\bfu^\top \bfS_{r_n}}\,
\1\big(a_n^{-1}M_{r_n}^{|\bfX|}\le x \big)\big]-1\big)\,(1+o(1))\,,\; 
\nto\,.
\eeam
By a telescoping sum argument we will prove that,
for every $\bfu\in\bbr^d$ and $x>0$,
\beam\label{eq:augu9}
\lefteqn{\lim_{k\to\infty}\limsup_{n\to\infty}
  \Big|
 k_n\,\big(
\E\big[\ex^{i\,a_n^{-1}\bfu^\top \bfS_{r_n}}\,
\1\big(a_n^{-1}M_{r_n}^{|\bfX|}\le x \big)\big]-1\big)}\nonumber\\
&&\hspace{0.7cm}-n\,\big( 
\E\big[\ex^{i\,a_n^{-1}\bfu^\top \bfS_{k}}\,\1(a_n^{-1}M_k^{|\bfX|}\le x)\big]
- \E\big[\ex^{i\,a_n^{-1}\bfu^\top \bfS_{k-1}}\,\1(a_n^{-1}M_{k-1}^{|\bfX|}\le x) \big] \big)\Big|=0\,.
\eeam
We write $\bfS_{-j}=\bfX_{-j}+\cdots +\bfX_{-1}$ and $M_{-j}^{|\bfX|}=\max_{-j\le t\le -1}|\bfX_t|$ 
for $j\ge 1$ as well as, for any integers $a\le b$, 
$M_{a,b}^{|\bfX|}=\max_{a\le t \le b}|\bfX_t|$.
The term inside the absolute values in \eqref{eq:augu9} can be expressed as the sum of
two negligible residual terms (which we omit in the sequel) 
and  the main term given by
\beam\label{eq:hot}
J_n:=
&&k_n\,\sum_{j=1}^{r_n-k}\E\big[\ex^{i\,a_n^{-1} \bfu^\top\bfS_{-k-j}}\1(a_n^{-1}M_{-k-j}^{|\bfX|}\le x)-\ex^{i\,a_n^{-1}\bfu^\top(
\bfS_{-k-j}- \bfS_{-j})}\1(a_n^{-1}M_{-k-j,-j-1}^{|\bfX|}\le x)\nonumber\\
&&\hspace{1.2cm}
-\ex^{i\,a_n^{-1}\bfu^\top \bfS_{-k-j+1}}\1(a_n^{-1}M_{-k-j+1}^{|\bfX|}\le x)\nonumber\\&&\hspace{1.2cm}+\ex^{i\,a_n^{-1}\bfu^\top( \bfS_{-k-j+1}- \bfS_{-j})}\1(a_n^{-1}M_{-k-j+1,-j-1}^{|\bfX|}\le x)\big]\,. 
\eeam 
Our goal is to approximate $J_n$ by some simpler quantities and to show 
that those are negligible as $\nto$ and $\kto$.
Suppressing the dependence on $x$, we write
\beao
A_n^c=\{a_n^{-1}|\bfX_{-k-j}|>x,a_n^{-1}M_{-j}^{|\bfX|}>x\}\,.
\eeao 
We first observe that
\beam
k_n\,\sum_{j=1}^{r_n-k} \P\big(A_n^c\big)
&=&k_n\,\sum_{j=1}^{r_n-k} \P\big(a_n^{-1}|\bfX_0|>x,a_n^{-1}M_{k,k+j-1}^{|\bfX|}>x\big)\nonumber\\
&\le &k_n\,\sum_{j=1}^{r_n-k} \sum_{l=k}^{k+j-1}\P\big(a_n^{-1}|\bfX_0|>x,a_n^{-1}|\bfX_l|>x\big)\nonumber\\
&\le & n\, \sum_{l=k}^{r_n}\P\big(a_n^{-1}|\bfX_0|>x,a_n^{-1}|\bfX_l|>x\big)\label{eq:augu1}\,.
\eeam
Following the discussion after the anti-clustering condition \eqref{cond:acac}, 
we conclude that the \rhs\ converges to zero by first letting
$\nto$ and then $\kto$. Hence, under the same limit regime, 
$J_n$ can be made arbitrarily close to
\beam
&&k_n\,\sum_{j=1}^{r_n-k}\E\Big[
\big(1-\1(A_n^c)\big)\label{eq:augu3}\\
&&\hspace{7mm}\times\Big(
\ex^{i\,a_n^{-1} \bfu^\top\bfS_{-k-j}}\1(a_n^{-1}M_{-k-j}^{|\bfX|}\le x)-\ex^{i\,a_n^{-1}\bfu^\top(
\bfS_{-k-j}- \bfS_{-j})}\1(a_n^{-1}M_{-k-j,-j-1}^{|\bfX|}\le x)\nonumber\\
&&\hspace{1cm}
-\ex^{i\,a_n^{-1}\bfu^\top \bfS_{-k-j+1}}\1(a_n^{-1}M_{-k-j+1}^{|\bfX|}\le x)\nonumber\\&&\hspace{1cm}+\ex^{i\,a_n^{-1}\bfu^\top( \bfS_{-k-j+1}- \bfS_{-j})}\1(a_n^{-1}M_{-k-j+1,-j-1}^{|\bfX|}\le x)\Big)\Big]\,.\nonumber
\eeam
%We observe that $1-\1(A_n^c)= \1(A_{n1})+\1(A_{n2})-\1(A_{n3})$ where
Writing
\beao
A_{n1}=\{a_n^{-1}|\bfX_{-k-j}|\le x\}\,,
A_{n2}=\{a_n^{-1}M_{-j}^{|\bfX|}\le x\}\,,
A_{n3}=
\{a_n^{-1}|\bfX_{-k-j}|\le x\,,a_n^{-1}M_{-j}^{|\bfX|}\le x\}\,,
\eeao
we may replace the multiplier $1-\1(A_n^c)$ in \eqref{eq:augu3}
by $\1(A_{n1})+\1(A_{n2})-\1(A_{n3})$. Then we obtain 
\beao
&&k_n\,\sum_{j=1}^{r_n-k}\Big(\E\big[\big(\ex^{i\,a_n^{-1} \bfu^\top\bfS_{-k-j}}-\ex^{i\,a_n^{-1}\bfu^\top \bfS_{-k-j+1}}\big)\1(a_n^{-1}M_{-k-j}^{|\bfX|}\le x)-\\
&&\hspace{1.4cm}
\big(\ex^{i\,a_n^{-1}\bfu^\top(
\bfS_{-k-j}- \bfS_{-j})}-\ex^{i\,a_n^{-1}\bfu^\top( \bfS_{-k-j+1}- \bfS_{-j})}\big)\1(a_n^{-1}M_{-k-j,-j-1}^{|\bfX|}\le x)\big]\\
&&\hspace{0.8cm}+\E\big[\big(\ex^{i\,a_n^{-1} \bfu^\top\bfS_{-k-j}}-\ex^{i\,a_n^{-1}\bfu^\top(
\bfS_{-k-j}- \bfS_{-j})}\big)\1(a_n^{-1}M_{-k-j}^{|\bfX|}\le x)-\\
&&
\hspace{1.4cm}\big(\ex^{i\,a_n^{-1}\bfu^\top \bfS_{-k-j+1}}-\ex^{i\,a_n^{-1}\bfu^\top( \bfS_{-k-j+1}- \bfS_{-j})}\big)\1(a_n^{-1}M_{-k-j+1}^{|\bfX|}\le x)\big]\\
&&\hspace{0.8cm}- \E\big[(1-\ex^{-i\,a_n^{-1}\,\bfu^\top S_{-j}})\,(\ex^{i\,a_n^{-1}\bfu^\top S_{-k-j}}
-\ex^{i\,a_n^{-1}\bfu^\top S_{-k-j+1}})\,\1(a_n^{-1}\,M_{-k-j}^{|\bfX|}\le x)
\big]
\Big)\\
&=:&I_{n1}+I_{n2}-I_{n3}\,, 
\eeao
and $I_{n1}$ turns into
 \beao
I_{n1}&=&k_n\,\sum_{j=1}^{r_n-k}\E\big[\big(\ex^{i\,a_n^{-1}\bfu^\top(
\bfS_{-k-j}- \bfS_{-j})}-\ex^{i\,a_n^{-1}\bfu^\top( \bfS_{-k-j+1}- \bfS_{-j})}\big)\\
&&\hspace{1.4cm}\big(\ex^{i\,a_n^{-1} \bfu^\top\bfS_{-j}}\,\1(a_n^{-1}M_{-j}^{|\bfX|}\le x)-1\big)\1(a_n^{-1}M_{-k-j,-j-1}^{|\bfX|}\le x)\big]\\
&=&k_n\,\sum_{j=1}^{r_n-k}\E\big[\big(\ex^{-i\,a_n^{-1}\bfu^\top \bfS_{-j}}
\big(\ex^{i\,a_n^{-1}\bfu^\top
\bfS_{-k-j}}- \ex^{i\,a_n^{-1}\bfu^\top\bfS_{-k-j+1}}\big)\big)\\
&&\hspace{1.4cm}\big((\ex^{i\,a_n^{-1} \bfu^\top\bfS_{-j}}-1)\,\1(a_n^{-1}M_{-j}^{|\bfX|}\le x) - \1(a_n^{-1}M_{-j}^{|\bfX|}>x)\big)\\&&\hspace{1.4cm}\1(a_n^{-1}M_{-k-j,-j-1}^{|\bfX|}\le x)\big]\,.
\eeao
Using   stationarity and the fact that $x\mapsto\ex^{i\,x}$ is a Lipschitz function 
bounded by 1, we obtain for $\bfu\ne\bf0$,
\beao
|I_{n1}|& \le &%k_n\,\sum_{j=1}^{r_n-k}\E\big[(|a_n^{-1}\bfu^\top \bfX_{-k-j}|\wedge
%2)\;\big|\ex^{i\,a_n^{-1} \bfu^\top\bfS_{\red -j}}\1(a_n^{-1}M_{-j}^{|\bfX|}\le x)-1\big|\big]\\
k_n\,\sum_{j=1}^{r_n-k} \E\Big[(|a_n^{-1}\bfu^\top \bfX_{-k-j}|\wedge
2)\,\Big((|a_n^{-1}\bfu^\top \bfS_{-j}|\wedge
2)+\1(a_n^{-1}M_{-j}^{|\bfX|}>x)\Big)\Big]\\
&\le&  n\, |\bfu|^2\sum_{l=k}^{r_n} \E\big[(|a_n^{-1}\bfX_{0}|\wedge
(3/|\bfu|))\,(|a_n^{-1} \bfX_{l}|\wedge
(3/|\bfu|))\big]\\&&+
n\,|\bfu|\,\sum_{l=k}^{r_n} \E\big[(|a_n^{-1}\bfX_{0}|\wedge
(3/|\bfu|))\,\1(a_n^{-1}|\bfX_l|>x)\big]\\
&\le &  n \,|\bfu|^2\sum_{l=k}^{r_n} \E\big[(|a_n^{-1}\bfX_{0}|\wedge
(3/|\bfu|))\,(|a_n^{-1} \bfX_{l}|\wedge
(3/|\bfu|))\big]\\
&&+ n\,|\bfu|\,\sum_{l=k}^{r_n} \E\big[(|a_n^{-1}\bfX_{0}|\wedge
(3/|\bfu|))\,(|(xa_n)^{-1}\bfX_l|\wedge 1\big)\big]\,.
\eeao
Here we also used  
sub-additivity.
Under the anti-clustering condition  \eqref{cond:acac} the \rhs\
converges to zero by first letting $\nto$ and then $\kto$.
\par 
We also have
\beao
I_{n2}&=&k_n\,\sum_{j=1}^{r_n-k}\E\big[\ex^{i\,a_n^{-1}\bfu^\top
\bfS_{-k-j+1}}\,\big(1-\ex^{-i\,a_n^{-1}\bfu^\top \bfS_{-j}}\big)\\
&&\hspace{1.2cm}\big(\big(\ex^{i\,a_n^{-1} \bfu^\top\bfX_{-k-j}}-1\big)\,\1(a_n^{-1}|\bfX_{-k-j}|\le x)-\1(a_n^{-1}|\bfX_{-k-j}|> x)\big)\\&&\hspace{1.2cm}\1(a_n^{-1}M_{-k-j+1}^{|\bfX|}\le x)\big]\,.
\eeao
Now, a similar argument as for $I_{n1}$ shows that $|I_{n2}|$ is negligible
by first letting $\nto$ and then $\kto$, and a similar argument also applies to
$|I_{n3}|$. 
Thus we proved \eqref{eq:augu9}. 
\par
Recalling \eqref{eq:augu10}--\eqref{eq:augu9}, it remains to characterize the limit $\Psi(\bfu,x)$:
\beao
&&\log\Psi(\bfu,x)\\
&=&\lim_{\kto}\lim_{\nto}
n\,\big( 
\E\big[\ex^{i\,a_n^{-1}\bfu^\top \bfS_{k}}\,\1(a_n^{-1}M_k^{|\bfX|}\le x)\big]
- \E\big[\ex^{i\,a_n^{-1}\bfu^\top \bfS_{k-1}}\,\1(a_n^{-1}M_{k-1}^{|\bfX|}\le x) \big] \big)\\
&=&
\lim_{\kto}\lim_{\nto} n\,
\E\big[\ex^{i\,a_n^{-1}\bfu^\top (\bfX_0+\bfS_k)}\,\1\big(a_n^{-1}(|\bfX_{0}|\vee M_{k}^{|\bfX|})\le x \big)\\&&\hspace{2.4cm}- \ex^{i\,a_n^{-1}\bfu^\top \bfS_k}\,\1\big(a_n^{-1}M_{k}^{|\bfX|}
\le x \big)\big]\nonumber\\
&=&\lim_{\kto}\lim_{\nto} n\,\big(\E\big[\ex^{i\,a_n^{-1}\bfu^\top (\bfX_0+\bfS_k)}\big]-\E\big[ \ex^{i\,a_n^{-1}\bfu^\top \bfS_k} \big]\big)\nonumber\\
&&+\lim_{\kto}\lim_{\nto} n\, \E\big[\ex^{i\,a_n^{-1}\bfu^\top \bfS_k}\,\1\big(a_n^{-1}M_{k}^{|\bfX|}> x \big)\nonumber\\
&&\hspace{3cm}-\ex^{i\,a_n^{-1}\bfu^\top (\bfX_0+\bfS_k)}\,\1\big(a_n^{-1}(|\bfX_{0}|\vee M_{k}^{|\bfX|})> x \big)\big]\nonumber\\
&=:& I_4+I_5\,,\qquad \bfu\in\bbr^d\,,\qquad x>0\,. 
\eeao
We will employ the \regvar\ and stationarity of $(\bfX_t)$ to deal with $I_4$ and $I_5$. In the last part of the proof we will show that
\beam\label{eq:I3}
I_4=\log\varphi_{\bfxi_\a}(\bfu)\,.
\eeam
As regards $I_5$, we observe that as $\nto$,
\beao
\lefteqn{n\, \E\big[\ex^{i\,a_n^{-1}\bfu^\top \bfS_k}\,
\1\big(a_n^{-1}M_{k}^{|\bfX|}> x \big)-\ex^{i\,a_n^{-1}\bfu^\top (\bfX_0+\bfS_k)}\,\1\big(a_n^{-1}(|\bfX_{0}|\vee M_{k}^{|\bfX|})> x \big)\big]}\nonumber\\
&=&\int_{\bbr^{dk}}\ex^{i\,\bfu^\top \sum_{t=0}^{k-1} \bfx_t}\,\1\Big(\max_{0\le t\le k-1}|\bfx_t|> x \Big)\, \big[n\,\P(a_n^{-1}(\bfX_0,\ldots,\bfX_{k-1})\in d\bfx )\big]\\&&-\int_{\bbr^{d(k+1)}} \ex^{i \bfu^\top \sum_{t=0}^{k} \bfx_t}\,\1\Big(\max_{0\le t\le k}|\bfx_t|> x \Big)\, \big[n\,\P(a_n^{-1}(\bfX_0,\ldots,\bfX_{k})\in d\bfx )\big]\\
&\to&\int_{\bbr^{dk}}\ex^{i\,\bfu^\top \sum_{t=0}^{k-1} \bfx_t}\,\1\Big(\max_{0\le t\le k-1}|\bfx_t|> x \Big)\, \mu_{k-1}(d\bfx )\\&&-\int_{\bbr^{d(k+1)}} \ex^{i \bfu^\top \sum_{t=0}^{k} \bfx_t}\,\1\Big(\max_{0\le t\le k}|\bfx_t|> x \Big)\, \mu_k(d\bfx )\\
&=&\int_{\bbr^{d(k+1)}}\Big[\ex^{i\,\bfu^\top \sum_{t=1}^{k} \bfx_t}\,\1\Big(\max_{1\le t\le k}|\bfx_t|> x \Big)-\ex^{i \bfu^\top \sum_{t=0}^{k} \bfx_t}\,\1\Big(\max_{0\le t\le k}|\bfx_t|> x \Big)\,\Big]\, \mu_k(d\bfx )
\,.
\eeao
In the limit relation we used the vague \con\ to the tail \ms\ 
$\mu_k(\cdot)$ in $\bbr^{d(k+1)}_{\bf0}$ defined in \eqref{eq:vagconv} 
for $k=0,1,\ldots$, and the facts that 
the integrands vanish in some neighborhood of the origin and are 
continuous with respect to any homogeneous measure.
The last identity follows from the consistency property 
of the tail \ms s $(\mu_k)$. Now, expressing the 
tail \ms s in terms of the tail process $(Y\,\bfTh_t)$ of $(\bfX_t)$ and using dominated \con , $I_5$ turns into
\beao
I_5&=& \lim_{k\to\infty}\int_0^\infty\E\Big[\ex^{i\,y\,\bfu^\top \sum_{t=1}^{k} \bfTh_t}\,\1\Big(y\,\max_{1\le t\le k}|\bfTh_t|> x \Big)\nonumber\\
&&\hspace{1.6cm}-\ex^{iy\,\bfu^\top \sum_{t=0}^{k} \bfTh_t}\,\1\Big(y\,\max_{0\le t\le k}|\bfTh_t|> x \Big)\Big]\,d\big(-y^{-\a}\big)\nonumber\\
&=&\int_0^\infty\E\Big[\ex^{iy\,\bfu^\top \sum_{t=1}^\infty \bfTh_t}\,\1\Big(y\,\max_{ t\ge 1}|\bfTh_t|>x \Big)\nonumber\\
&&\hspace{1cm}-\ex^{i\,y\bfu^\top \sum_{t=0}^\infty \bfTh_t}\,\1\Big(y\,\max_{t\ge 0}|\bfTh_t|> x \Big)\Big]\,d\big(-y^{-\a}\big)\,.
\eeao
Using the definition of $\bfQ=\bfTh/\|\bfTh\|_\a$ together with the time-change formula \eqref{eq:tcp}, the \rhs\ can be written as 
\beao
%\lefteqn{\int_0^\infty\Big(\E\Big[\ex^{iy\,\bfu^\top \sum_{t=0}^\infty \bfTh_t}\,\1\Big(y\,\max_{ t\ge 0}|\bfTh_t|>x \Big)}\nonumber\\
%&&\hspace{1cm}-\ex^{i\,y\bfu^\top \sum_{t= 1}^\infty \bfTh_t}\,\1\Big(y\,\max_{t\ge 1}|\bfTh_t|>x \Big)\Big]\Big)\,d\big(-y^{-\a}\big)\nonumber\\
I_5&=&\int_0^\infty\sum_{j\in\bbz}\E\Big[|\bfTh_j|^\a\Big(\ex^{iy\,\bfu^\top 
\sum_{t=1}^\infty \bfQ_t}\,\1\Big(y\,\max_{t\ge 1}|\bfQ_t|> x \Big)\nonumber\\
&&\hspace{2.3cm}-\ex^{i\,y\bfu^\top \sum_{t=0}^\infty \bfQ_t}\,\1\Big(y\,\max_{t\ge 0}|\bfQ_t|> x \Big)\Big)\Big]\,d\big(-y^{-\a}\big)\nonumber\\
&=&\int_0^\infty\sum_{j\in\bbz}\E\Big[ \Big(\ex^{iy\,\bfu^\top \sum_{t=-j+1}^\infty \bfQ_t}\,\1\Big(y\,\max_{ t\ge -j+1}|\bfQ_t|> x \Big)\nonumber\\
&&\hspace{2.3cm}-\ex^{i\,y\bfu^\top \sum_{t=-j}^\infty \bfQ_t}\,\1\Big(y\,\max_{t\ge -j}|\bfQ_t|> x \Big)\Big)\Big]\,d\big(-y^{-\a}\big)\nonumber\\
&=&-\dint_0^\infty\E\Big[\ex^{i\,y\,\bfu^\top \sum_{t=-\infty}^\infty \bfQ_t}\, 
\1\Big(y\,\max_{t\in \bbz}|\bfQ_t|> x \Big)\Big]\,d(-y^{-\a})\,.
\eeao
Thus we proved that 
\beao
\log\Psi(\bfu,x)&=&I_4+I_5\\
&=&\log\varphi_{\bfxi_\a}(\bfu)-\dint_0^\infty\E\Big[\ex^{i\,y\,\bfu^\top \sum_{t=-\infty}^\infty \bfQ_t}\, 
\1\Big(y\,\max_{t\in \bbz}|\bfQ_t|> x \Big)\Big]\,d(-y^{-\a})\,.
\eeao
This is the logarithm of the desired limit hybrid \chf\ in terms of $\bfQ$.
Changing variables and appealing to the definition of $\wt \bfQ$
in  \eqref{eq:qtilde}, we also have
\beao
&&\Psi(\bfu,x)\\
&=&\varphi_{\bfxi_\a}(u)\exp\Big(-\dint_0^\infty\Big(\E\Big[\Big(\ex^{i\,y\,\bfu^\top\sum_{t=-\infty}^\infty \bfQ_t}-1\Big)\,\1\Big(y\,\max_{t\in \bbz}|\bfQ_t|>x \Big)\Big]\, d(-y^{-\a})\Big)\\
&& \times\exp\Big( - \dint_0^\infty\P\Big(y\,\max_{t\in \bbz}|\bfQ_t|>x \Big) \,d(-y^{-\a})\Big)\\
&=& \varphi_{\bfxi_\a}(u)\,\Phi_\a^{\theta_{|\bfX|}}(x)\\
&&\times
\exp\Big(-\theta_{|\bfX|}\dint_x^\infty\Big(\E\Big[\dfrac{\max_{t\in \bbz}|\bfQ_t|^\a}
{\theta_{|\bfX|}}\,\Big(\ex^{i\,y\,\bfu^\top\sum_{t=-\infty}^\infty 
\bfQ_t/\max_{t\in \bbz}|\bfQ_t|}-1\Big)\Big]\,d(-y^{-\a})\Big) 
\\
&=& \varphi_{\bfxi_\a}(u)\,\Phi_\a^{\theta_{|\bfX|}}(x)\, \exp\Big(-\theta_{|\bfX|}\dint_x^\infty\Big(\E\Big[ \ex^{i\,y\,\bfu^\top\sum_{t=-\infty}^\infty \wt \bfQ_t}-1 \Big]\,d(-y^{-\a})\Big)\,, 
\eeao
and the desired result follows.\\[2mm]
{\em Proof of \eqref{eq:I3}.}
An application of the \cmt\ for \regvary\ random vectors
ensures that $\bfS_k$ inherits \regvar\ with index $\a$:
for $k\ge 1$,
\beao
n\,\P(a_n^{-1}\bfS_k\in\cdot )
&\stv& \wt  \mu_{\bfS_k}(\cdot) := \mu_{k-1}\big(\{(\bfx_0,\ldots,\bfx_{k-1})\in \R^{dk}:\, \bfx_0+\cdots+\bfx_{k-1}\in\cdot \}\big)\,,\qquad \nto\,.
\eeao 
Therefore the \ds\ of $\bfS_k$  belongs to the domain of attraction
of an $\alpha$-stable  law (possibly degenerate), hence there exists an
$\alpha$-stable \rv\ $\bfxi_\alpha^{(k)}$ \st\ for iid copies $(\bfS_{k,i})$ of $\bfS_k$,
\beao 
(a_n^{(k)})^{-1}\sum_{i=1}^n \bfS_{k,i}\std \bfxi_\alpha^{(k)}\,,
\eeao
where 
\beao
a_n^{(k)}=a_n\,\big[\mu_{k-1}\big(\{(\bfx_0,\ldots,\bfx_{k-1})\in \R^{dk}:\, |\bfx_0+\cdots+\bfx_{k-1}|>1 \}\big)\big]^{1/\a}=:a_n\,C_k^{1/\a}\,,\eeao 
\st
\beao
n\,\P(|\bfS_k|>a_n^{(k)})
&\sim& n\,\P(|\bfX|>a_n^{(k)})\,C_k\to 1\,,\qquad\nto\,.
\eeao
Since $\bfS_k$ is \regvary\ one can also define its tail measure as the vague limit 
 \beao
 \dfrac{\P(x^{-1}\bfS_k\in\cdot )}{\P(|\bfS_k|>x)}
= \dfrac{\P(|\bfX|>x)}{{\P(|\bfS_k|>x)}}\dfrac{\P(x^{-1}\bfS_k\in\cdot )}{\P(|\bfX|>x)}
\stv  \mu_{\bf S_k}(\cdot)=C_k^{-1} \wt\mu_{\bfS_k}(\cdot)\,.
\eeao
An application of  Lemma 3.5 in Petrov
\cite{petrov:1995} yields the equivalent relation
\beam
n\,\big(\varphi_{a_n^{-1} \bfS_{k} }(\bfu)-1\big) =n\,\big(\varphi_{(a_n^{(k)}/a_n)\,(a_n^{(k)})^{-1} \bfS_{k} }(\bfu)-1\big)\to \log  \varphi_{C_k^{1/\a}\bfxi_\alpha^{(k)}}(\bfu)\,,\qquad \bfu\in\bbr^d\,.\nonumber\\
\label{eq:xmas30b}
\eeam
The log-\chf\ of the $\a$-stable $\bfxi_\alpha^{(k)}$ 
can be written as 
\beao
\log  \varphi_{ C_k^{1/\a}\bfxi_\alpha^{  (k)}}(\bfu) &=&\log  \varphi_{\bfxi_\alpha^{(k)}}(C_k^{1/\a}\, \bfu) \\
&=&  \dint_{\bbr_{\bf0} ^d}\big( \ex^{i\,C_k^{1/\a}\,\bfu^\top\, \bfy} -1 - i\,C_k^{1/\a}\,\bfu^\top\,\bfy\;  \1_{(1,2)}(\alpha)\big)\,\mu_{\bfS_k}(d\bfy)\\
&=:& -c_\alpha\, \sigma_k^\a(C_k^{1/\a}\,\bfu)\,\big(1
-i\,\beta_k(C_k^{1/\a}\,\bfu) \,\tan(\a\,\pi/2)\big)\,,
\qquad \bfu\in\bbr^d\,,
\eeao
where, by homogeneity of the tail measure,
\beao
\sigma_k^\a(C_k^{1/\a}\,\bfu)&:=&\mu_{\bfS_k}\big(\big\{\bfx\in \R^d:\, |{ C_k^{1/\a}\, \bfu^\top \bfx}| >1\big\}\big)=\wt \mu_{\bfS_k}\big(\big\{\bfx\in \R^d:\, | \bfu^\top \bfx| >1\big\}\big)\,,\\
\beta_k(C_k^{1/\a}\bfu)&:=&\dfrac{\wt \mu_{\bfS_k}\big(\big\{\bfx\in \R^d:\,  \bfu^\top \bfx  >1\big\}\big)-\wt \mu_{\bfS_k}\big(\big\{\bfx\in \R^d:\,  \bfu^\top \bfx  <-1\big\}\big)}
{\wt \mu_{\bfS_k}\big(\big\{\bfx\in \R^d:\, | \bfu^\top \bfx| >1\big\}\big)}\,.
\eeao
Then, using the definition of $\wt \mu_{\bfS_k}$ and the consistency of the tail \ms s $(\mu_k)$, we have
\beao
\sigma_k^\a(C_k^{1/\a}\bfu)&=&
%\mu_{\bfS_k}\Big(\Big\{\bfx\in \R^d;\, | \bfu^\top \bfx| >1\Big\}
\mu_{k-1} \big(\{(\bfx_0,\ldots,\bfx_{k-1})\in \bbr^{d k}: 
|\bfu^\top (\bfx_0+\cdots+\bfx_{k-1})|>1 \}\big)\\
&=&\mu_{k} \big(\{(\bfx_0,\ldots,\bfx_{k})\in \bbr^{d  (k+1)}: 
|\bfu^\top (\bfx_0+\cdots+\bfx_{k-1})|>1 \}\,,\bfx_k\in\bbr^d\}\big)\\
&=&\mu_{k}\big(\{(\bfx_0,\ldots,\bfx_{k})\in \bbr^{d (k+1)}: 
|\bfu^\top (\bfx_1+\cdots+\bfx_{k})|>1 \}\,,\bfx_0\in\bbr^d\}\big)\,,
\eeao
and the quantities $\beta_k(\bfu)$ can be expressed similarly. 
Keeping this remark in mind and exploiting \eqref{eq:xmas30b}, we have as $\nto$,
\beao
I_4&=&
n\, \big(\varphi_{a_n^{-1} \bfS_{k+1} }(\bfu) -\varphi_{a_n^{-1} \bfS_{k} }(\bfu) \big)\nonumber\\[2mm]
&\to& \log \varphi_{C_{k+1}^{1/\a}\,\bfxi_\alpha^{(k+1)}}(\bfu)-\log \varphi_{C_{k}^{1/\a}\,\bfxi_\alpha^{(k)}}(\bfu)\\[2mm]
&=&   -c_\alpha \,\Big(\underbrace{\big(\sigma_{k+1}^\a(C_{ k+1}^{1/\a}\bfu)-\sigma_k^\a(C_k^{1/\a}\bfu)\big)}_{=:\Delta_1(k)}\\
%\Big( \int_{\R^{dk}_{\bf0}}\1\Big(  \Big| \bfu^\top \sum_{i=0}^k\bfx_i\Big| >1 \Big)-\1\Big(  \Big|\bfu^\top \sum_{i=1}^{k}\bfx_i\Big| >1 \Big)\mu_{k}(d\bfx )\label{eq:coeff1}\\
&&\hspace{1cm}-i\,\tan( \a\,\pi/2) \underbrace{\big(\beta_{k+1}(C_{ k+1}^{1/\a}\bfu)\sigma_{k+1}^\a(C_{ k+1}^{1/\a}\bfu)- \beta_{k}(C_k^{1/\a}\bfu)\sigma_{k}^\a(C_k^{1/\a}\bfu)\big)}_{=:\Delta_2(k)}\Big)\,.
%\Big( \int_{\R^{dk}_{\bf0}}\1\Big(    \bfu^\top \sum_{i=0}^k\bfx_i  >1 \Big)-\1\Big(   \bfu^\top \sum_{i=1}^{k}\bfx_i >1 \Big)\mu_{k}(d\bfx ) \label{eq:coeff2}\\
%&&-\int_{\R^{dk}_{\bf0}}\1\Big(    \bfu^\top \sum_{i=0}^k\bfx_i  <-1 \Big)-\1\Big(   \bfu^\top \sum_{i=1}^{k}\bfx_i <-1 \Big)\mu_{k}(d\bfx )\Big)\Big)\label{eq:coeff3}\,.
\eeao
Our next goal is to identify $\Delta_1(k)$ and $\Delta_2(k)$ as follows:
\beao
\Delta_1(k)&=&\E\Big[  \Big| \bfu^\top \sum_{i=0}^k\bfTh_i \Big|^\a- \Big| \bfu^\top \sum_{i=1}^{k} \bfTh_i\Big|^\a\Big]\,,\\
\Delta_2(k)&=&\E\Big[ \Big( \Big( \bfu^\top \sum_{i=0}^k\bfTh_i \Big)_+^\a-\Big( \bfu^\top \sum_{i=1}^k\bfTh_i \Big)_+^\a\Big)-\Big(\Big( \bfu^\top \sum_{i=0}^{k} \bfTh_i\Big)_-^\a- \Big( \bfu^\top \sum_{i=1}^{k} \bfTh_i\Big)_-^\a\Big)\Big]\,.
\eeao
We give the details only for $\Delta_1(k)$.
Rewriting the tail \ms s $\mu_k$ in terms of the tail process $(Y\,\bfTh_t)$, we have
\beao
\Delta_1(k)&=&
%\lefteqn{\int_{\R^{dk}_{\bf0}}\1\Big(    \Big|\bfu^\top \sum_{i=0}^k\bfx_i  \Big|>1 \Big)-\1\Big(  \Big| \bfu^\top \sum_{i=1}^{k}\bfx_i\Big| >1 \Big)\mu_{k}(d\bfx )}\\
\int_{\R^{dk}_{\bf0}}\Big(\1\Big(   \Big| \bfu^\top \sum_{i=0}^k\bfx_i\Big|  >1 \Big)-\1\Big(  \Big| \bfu^\top \sum_{i=1}^{k}\bfx_i \Big|>1 \Big)\Big)\,\1(|\bfx_0|>0)\,\mu_k(d\bfx)\\
&=&\int_0^\infty\Big(\P\Big( y \,  \Big|\bfu^\top \sum_{i=0}^k\bfTh_i\Big|  >1 \Big)-\P\Big(  y\,\Big| \bfu^\top \sum_{i=1}^{k}\bfTh_i \Big|>1 \Big)\Big)\,d(-y^{-\a})\\
&=&\E\Big[  \Big| \bfu^\top \sum_{i=0}^k\bfTh_i \Big|^\a- \Big| \bfu^\top \sum_{i=1}^{k} \bfTh_i\Big|^\a\Big]\,.
\eeao
\par
Our next step is to show that the following limits exist
\beam
\Delta_1:=\lim_{\kto} \Delta_1(k)&=&\E\Big[   \Big| \bfu^\top \sum_{i=0}^\infty\bfTh_i \Big|^\a- \Big| \bfu^\top \sum_{i=1}^{\infty} \bfTh_i \Big|^\a\Big]\label{eq:august21c}
\\
\Delta_2:=\lim_{\kto}\Delta_2(k)&=& 
\E\Big[ \Big(  \Big( \bfu^\top \sum_{i=0}^\infty\bfTh_i \Big)_+^\a- \Big( \bfu^\top \sum_{i=1}^{\infty} \bfTh_i \Big)_+^\a\Big)\nonumber\\&&  \;\;- \Big(\Big( \bfu^\top \sum_{i=0}^\infty\bfTh_i \Big)_-^\a- \Big( \bfu^\top \sum_{i=1}^{\infty} \bfTh_i \Big)_-^\a\Big)
\Big]\,.\nonumber
\eeam
Then  the limit
\beao
\log \varphi_{C_{k+1}^{1/\a}\,\bfxi_\alpha^{(k+1)}}(\bfu)-\log \varphi_{C_{k}^{1/\a}\,\bfxi_\alpha^{(k)}}(\bfu)\to \log \varphi_{\bfxi_\alpha}(\bfu)\,,\qquad\kto\,,
\eeao
exists and is equal to $  -c_\alpha(\Delta_1 -i\,\tan( \a\,\pi/2)  \Delta_2)$. This expression agrees with the \chf\ $\varphi_{\bfxi_\a}$ given in \eqref{eq:agree}, \eqref{eq:agree1}.
We restrict ourselves to the calculations for $\Delta_1$; the case
$\Delta_2$ is similar.
Using the following con\seq\ of the  mean value theorem  
\beao
|a+b|^\a-|b|^\a \le (\alpha\vee 1)(|a|+|b|)^{(\a-1)\vee 0}|a|^{\a\wedge 1}
\eeao and the fact that $|\bfTh_0|=1$, we have  for fixed $\bfu\in\bbr^d$,
\beam\label{eq:xmas30c}
  |\Delta_1(k)|
%\E\Big[   \Big| \bfu^\top \sum_{i=0}^k\bfTh_i \Big|^\a- \Big| \bfu^\top \sum_{i=1}^{k} \bfTh_i \Big|^\a\Big]
&\le&c\;\E\Big[  \Big(1+ \Big| \sum_{i=1}^k \bfu^\top \bfTh_i \Big|\Big)^{(\a-1)\vee 0}\Big]\nonumber\\
&\le&c\; \E\Big[ \Big( \sum_{i=0}^\infty |  \bfTh_i |\Big)^{(\a-1)\vee0}\Big] <\infty\,.
\eeam
%In the last step we used Jensen's inequality for $\a\in (1,2)$, the case
%$\alpha\in (0,1)$ being trivial. 
The \rhs\ is trivially finite for $\a\in (0,1)$ while for $\a\in (1,2)$
we employ the anti-clustering condition \eqref{cond:acac}; see Lemma \ref{lem:xmas29a}.
In view of \eqref{eq:xmas30c} we are allowed to use
dominated convergence and conclude that \eqref{eq:august21c} holds.
%\beao
%\lim_{\kto}\Delta_1(k)= %\E\Big[   \Big| \bfu^\top \sum_{i=0}^k\bfTh_i \Big|^\a-  \Big| \bfu^\top \sum_{i=1}^{k} \bfTh_i \Big|^\a\Big]\std 
%\E\Big[   \Big| \bfu^\top \sum_{i=0}^\infty\bfTh_i \Big|^\a- \Big| \bfu^\top \sum_{i=1}^{\infty} \bfTh_i \Big|^\a\Big]=\Delta_1 \,.
%\eeao
In the limit 
we introduce the spectral cluster process $\bfQ=\bfTh/\|\bfTh\|_\a$:
\beao
\Delta_1&=&
%\E\Big[   \Big| \bfu^\top \sum_{i=0}^\infty\bfTh_i \Big|^\a- \Big| \bfu^\top \sum_{i=1}^{\infty} \bfTh_i \Big|^\a\Big]&=&
\E\Big[ \|\bfTh\|_\a^\a\, \Big( \Big| \bfu^\top \sum_{i=0}^\infty\dfrac{\bfTh_i}{\|\bfTh\|_\a} \Big|^\a- \Big| \bfu^\top \sum_{i=1}^{\infty} \dfrac{\bfTh_i}{\|\bfTh\|_\a} \Big|^\a\Big)\Big]\\
&=&\sum_{t\in\bbz}\E\Big[ |\bfTh_t|^\a \Big( \Big| \bfu^\top \sum_{i=0}^\infty\dfrac{\bfTh_i}{\|\bfTh\|_\a} \Big|^\a- \Big| \bfu^\top \sum_{i=1}^{\infty} \dfrac{\bfTh_i}{\|\bfTh\|_\a} \Big|^\a\Big)\Big]\\
&=&\sum_{t\in\bbz}\E\Big[   \Big| \bfu^\top \sum_{i=-t}^\infty\dfrac{\bfTh_i}{\|\bfTh\|_\a} \Big|^\a- \Big| \bfu^\top \sum_{i=1-t}^{\infty} \dfrac{\bfTh_i}{\|\bfTh\|_\a} \Big|^\a\Big]\,.
\eeao
 The last identity follows by multiple application of the time-change 
formula \eqref{eq:tcp}.
We observe that we have a telescoping sum structure, resulting in 
\beao
\Delta_1&=&\E\Big[   \Big| \bfu^\top \sum_{i\in\bbz}\dfrac{\bfTh_i}{\|\bfTh\|_\a} \Big|^\a\Big]
=\E\Big[   \Big| \bfu^\top \sum_{i\in\bbz} \bfQ_i  \Big|^\a\Big]\,.
\eeao
A similar argument yields the identity
\beam \label{eq:xmasjan4a}
\E\Big[\Big(   \sum_{i\in\bbz} |\bfQ_i|  \Big)^\a\Big]=\E\Big[   \Big( \sum_{i=0}^\infty|\bfTh_i| \Big)^\a- \Big( \sum_{i=1}^{\infty} |\bfTh_i| \Big)^\a\Big]\,.
\eeam
In particular, the expectation on the \lhs\ is finite.  
\end{proof}

\subsection{Joint convergence of sums, maxima and $\ell^p$-norms}
Our next goal is to prove joint \con \ of $a_n^{-1}(\bfS_n,M_n^{|\bfX|},\gamma_{n,p})$. We work under the following mixing condition slightly stronger than \eqref{eq:wdepms}:

For some integer \seq s $r_n\to \infty$, $k_n=[n/r_n]\to\infty$,
\beam\label{eq:wdepms2}
\Psi_{n,p}(\bfu,x,\la)&:=&\E\Big[\exp\big(i\,a_n^{-1}\bfu^\top\bfS_{n}-a_n^{-p}\la \gamma_{n,p}^p\big)\,\1\big(a_n^{-1}M_{n}^{|\bfX|}\le x\big)\Big]\nonumber\\
 &=&
\Big(\E\Big[\exp\big(i\,a_n^{-1}\bfu^\top\bfS_{r_n}-a_n^{-p}\la \gamma_{r_n,p}^p\big)\,\1\big(a_n^{-1}M_{r_n}^{|\bfX|}\le x\big)\Big]\Big)^{k_n}+o(1),\nonumber\\&&\qquad \nto,\qquad (\bfu,x,\la)\in\bbr^d\times \R_+^2\,.
\eeam
This mixing condition can be derived by strong mixing or coupling arguments; see Section~\ref{sec:examples} below.
\bth\label{pr:prophybridch}
Assume the conditions of Theorem~\ref{th:joint} and the mixing condition \eqref{eq:wdepms2}. Then, with the notation of the latter result,  for $\a<p$,
\beao
a_n^{-1} (\bfS_n, M_n^{|\bfX|},\gamma_{n,p})\std (\bfxi_\a,\eta_\a,\zeta_{\a,p})\,,\qquad \nto\,,
\eeao
where the joint limit \ds\  
is described by the joint hybrid \chf --Laplace transform
of  $(\bfxi_\a,\eta_\a,\zeta_{\a,p}^p)$ given by
\beam\label{eq:xxx}\lefteqn{
\E\big[\ex^{i \,\bfu^\top\,\bfxi_\a}\,\1(\eta_\a\le x)\,\ex^{-\la\,\zeta_{\a,p}^p}\big]}\\\nonumber
&=&\exp\Big(\dint_0^\infty\E\Big[\ex^{i\,y\, \bfu^\top\sum_{t=-\infty}^\infty  \bfQ_t-y^p\la\sum_{t=-\infty}^\infty  |\bfQ_t|^p}\, 
\1\Big(y\,\max_{t\in \Z}|\bfQ_t|\le x \Big)\\
&&\nonumber-1-
i\,y\,\bfu^\top \sum_{t\in\Z} \bfQ_t\,\1_{(1,2)}(\a)\Big]\,d(-y^{-\a})\Big)\,,\qquad \bfu\in\R^d,\,x>0,\,\lambda>0\,.
%&=&\varphi_{\bfxi_\a}(\bfu)\,
%\exp\Big(-\dint_0^\infty\E\Big[\ex^{i\,y\,\bfu^\top \sum_{t=-\infty}^\infty Q_t}\, \1\Big(y\,\max_{t\in \Z}|Q_t|>x \Big)\Big]\,d(-y^{-\a})\Big)\nonumber\\
%&=&\varphi_{\bfxi_\a}(\bfu)\,\Phi_\a^{\theta_{|\bfX|}}(x)\, \exp\Big(-\theta_{|\bfX|}\dint_x^\infty\E\big[\ex^{i\,y\,\bfu^\top\sum_{t=-\infty}^\infty \wt \bfQ_t}-1\big]\,d(-y^{-\a})\Big)\nonumber\\
%&=&\varphi_{\bfxi_\a}(\bfu)\,\Phi_\a^{\theta_{|\bfX|}}(x)\,
%\exp\big(-\E\big[\ex^{i\,\bfu^\top Y\,\sum_{t=-\infty}^\infty \wt \bfQ_t}-1\big]\big)\,,\nonumber\\
\eeam
\ethe

Analogous results about the joint convergence of maxima, partial sums and finitely many $\gamma_{n,p_k}$ for $p_k>0 $ can be proved by similar techniques as below.
In Theorem~\ref{th:joint} the joint \con\ of sums and maxima 
was described via the \con\ of hybrid \chf s. Here we appeal to the 
\con\ of an extension of these: we use the point-wise \con\
of the joint hybrid \chf --Laplace transform of $a_n^{-1}(\bfS_n,M_n^{|\bfX|})$ and 
$a_n^{-p}\gamma_{n,p}^p$. 
\bre\label{rem:expnormfinite}\label{rem:stableskew}
Notice that $\E[\|\bfQ\|_p^\a]<\infty$ for $\a/p<1$ since, by concavity of the function $x\mapsto x^{\a/p}$, $x>0$,   we obtain $\|\bfQ\|_p^\a\le \|\bfQ\|_\a^\a=1$ a.s.
Therefore the limit Laplace transform of $(a_n^{-p}\gamma_{n,p}^p)$ is well defined.
We consider a Fr\'echet \rv\ $Y_p$ with index $p>0$
which is  independent of $\bfQ$.
Then the distribution of $\zeta_{\a,p}$ is provided by the Laplace transform
\beao
\E\big[\ex^{-\la\,\zeta_{\a,p}^p}\big]
&=&\exp\Big(\dint_0^\infty\E\Big[\ex^{-y^p\la\sum_{t=-\infty}^\infty  |\bfQ_t|^p} 
-1\Big]\,d(-y^{-\a})\Big)\\
&=&\exp\Big(-\dint_0^\infty\E\big[\P\big(y\,\la^{1/p}\,Y_p\,\|\bfQ\|_p>1
\mid \bfQ\big)\big]\,d(-y^{-\a})\Big)\\
&=&\exp\Big(-\E\Big[\dint_0^\infty\P\big(\la^{\a/p}\,Y_p^\a\,\|\bfQ\|_p^\a>u
\mid \bfQ\big)du\Big]\Big)\\
&=&\exp\big( -\E\big[Y_p^\a]\E\big[\|\bfQ\|_p^\a\big]\la^{\a/p} \big)\,.
\eeao
Recall the Laplace transform of a totally skewed to the right 
$\alpha/p$-stable  
distribution: 
\beao
\E\big[\ex^{-\la\,\xi_{\a/p}}\big]&=&\exp\big(-\Gamma(1-\alpha/p)\,\la^{\a/p}\big)\,,\qquad \la>0\,,
\eeao
see Samorodnitsky and Taqqu \cite{samorodnitsky:taqqu:1994}, 
Proposition~1.2.11.
Computing  $\E[Y_p^\a]=\Gamma(1-\alpha/p)$, we finally
get
\beao
 \E\big[\ex^{-\la\,\zeta_{\a,p}^p}\big]&=&\big( \E\big[\ex^{-\la\,\xi_{\a/p}}\big]\big)^{\E[\|\bfQ\|_p^\a]}\,,\qquad \la>0\,.
\eeao
Due to the presence of $\E[\|\bfQ\|_p^\a]\le 1$ in the latter expression we obtain a quantitative description of the clustering effect  in the dependent \seq\ $(\bfX_t)$ which is similar to the interpretation of the extremal index $\theta_{|\bfX|}$.
In particular, we observe that $ \E[\|\bfQ\|_p^\a]= 1$ in the \asy ally independent
case when $\bfQ_t=\bf0$ for $t\ne 0$.
\ere

\bre
The joint characteristic function\ -- Laplace transform $\Phi_{\bfxi_\a,\zeta_{\a,p}^p}(\bfu,\lambda)=\E\big[\ex^{i \bfu^\top \xi_\alpha-\lambda \zeta_{\a,p}^p}\big]$ 
is obtained by letting $x\to\infty$ in \eqref{eq:xxx}.
%$\E\big[\ex^{i \,\bfu^\top\,\bfxi_\a}\,\1(\eta_\a\le x)\,\ex^{-\la\,\zeta_{\a,p}^p}\big]$. 
The quantity $\Phi_{\bfxi_\a,\zeta_{\a,p}^p}(\bfu,\la)$ indicates that  $(\bfxi_\a,\zeta_{\a,p}^p)$ is not stable but it is still  infinitely divisible; this is a direct con\seq\ of the mixing condition and can be seen directly in the structure of $\Psi(\bfu,\lambda)$: 
for every $c>0$,
\beao
 \Phi_c(\bfu,\lambda):=\E\big[\ex^{i c\bfu^\top \bfxi_\alpha-\lambda c^{p}\zeta_{\a,p}^p}\big]
= (\Phi_{\bfxi_\a,\zeta_{\a,p}^p}(\bfu,\la))^{c^\a}\,.
%\exp\Big( c^\a \dint_0^\infty\E\Big[\ex^{i y \bfu^\top \sum_{t=-\infty}^\infty \bfQ_t -y^p\la\sum_{t=-\infty}^\infty  |\bfQ_t|^p} 
%-1\Big]\,d(-y^{-\a})\Big) 
\eeao
Thus, for every integer $n\ge 1$, $\Phi_{\bfxi_\a,\zeta_{\a,p}^p}(\bfu,\lambda)=(\Phi_{n^{-1/\a}}(\bfu,\la))^n$. This divisibility is  also observed in the joint characteristic function\ -- Laplace transform of norms $\E\big[\ex^{i\,u\,\bfxi_{\a/q}-\la\,\zeta_{\a/q,r/q}^{r/q}}\big]$; see
Section~\ref{subsec:ration:norms} below. Moreover, $\Phi_{\bfxi_\a,\zeta_{\a,p}^p}(\bfu,\lambda)$ does not describe a self-decomposable distribution. However, the following self-decomposition is valid:
\[
 \Phi_{\bfxi_\a,\zeta_{\a,p}^p}(\bfu,\lambda)=\Phi_{\bfxi_\a,\zeta_{\a,p}^p}(c\bfu,c^p \la)\,\Phi_{(1-c^{\a})^{-1/\a}}(\bfu,\la)\quad \text{for}\quad c\in(0,1). 
\]
\ere

\begin{proof}
We consider the joint hybrid \chf\ -- Laplace transform 
of $a_n^{-1}(\bfS_{r_n},M_{r_n}^{|\bfX|})$  and $(a_n^{-p}\gamma_{r_n,p}^p)$
given by
\beao
\Psi_{n,p}(\bfu,\la,x)=\E\big[\exp\big(i\,a_n^{-1}\bfu^\top\bfS_{r_n}-a_n^{-p}\la \gamma_{r_n,p}^p\big)\,\1\big(a_n^{-1}M_{r_n}^{|\bfX|}\le x\big)\big]\,.
\eeao
Consider iid copies 
$(Y_{p,t})$  of a Fr\'echet \rv\ $Y_p$ with \ds\ $\Phi_p$ which are also independent of $(\bfX_t)$ and construct the 
partial maxima
\beao
\wt M_{r_n}= \max_{1\le t\le r_n} Y_{p,t}|\bfX_t| \,,\qquad n\ge 1\,.
\eeao
Then we obtain for
$(\bfu,\la,x)\in\bbr^d\times\bbr_+^2$,
\beao
\lefteqn{\E\big[\ex^{i \,a_n^{-1}\,\bfu^\top\,\bfS_{r_n}\,
-\la\, a_n^{-p} \gamma_{r_n,p}^p} \,\1\big(a_n^{-1}M_{r_n}^{|\bfX|}\le x\big)\big]}\\
%&=&\E\big[\ex^{i \,a_n^{-1}\,u\,S_n}\,\1(a_n^{-1}M_n^{|X|}\le x)\,\ex^{i\,\la \,Y_p\,
%a_n^{-1} \gamma_{n,p}}\big]\\
%&=&\E\big[\ex^{i \,a_n^{-1}\,\bfu^\top\,\bfS_n}\,\1(a_n^{-1}M_n^{|\bfX|}\le x)\,\ex^{i\,\la \,Y_{p,1}\,
%a_n^{-1} \gamma_{n,p}}\big]\\
&=&\E\big[\ex^{i \,a_n^{-1}\,\bfu^\top\,\bfS_{r_n}}\,\1\big(a_n^{-1}\wt M_{r_n}\le \la^{-1/p}\,, a_n^{-1}M_{r_n}^{|\bfX|}\le x\big)\big]\,.
\eeao
We interpret the \rhs\ as hybrid \chf\  of the process
\beao
\bfZ_t:=(\bfX_t, Y_{p,t}|\bfX_t|,|\bfX_t|)\,,\qquad t\in\bbz\,.
\eeao
Since
$p>\a$ the \rv\ $Y_{p,t}$ has moments of order $p'\in (\a,p)$. Therefore an 
application of the multivariate 
Breiman theorem of Basrak et al. \cite{basrak:davis:mikosch:2002}
and \regvar\ of $(\bfX_t)$ 
yield that 
$(\bfZ_t)$ is an $\bbr^{d+2}$-valued \regvary\ 
stationary \seq\ with index $\a$. Writing $(\bfTh_t)$ for the spectral
tail process of $(\bfX_t)$ and exploiting the independence of 
$(\bfX_t)$ and $(Y_{p,t})$, we also have for $h\ge 0$,
\beao\lefteqn{\left.
\P\left(\dfrac 1 {|\bfX_0|} 
\left(
\left(\barr{rr} \bfX_0\\Y_{p,0}|\bfX_0|\\|\bfX_0|\earr\right)\,,\ldots,
\left(\barr{rr}  \bfX_h\\Y_{p,h}|\bfX_h|\\|\bfX_h|\earr\right)
\right)\in\cdot\,\right|\,|\bfX_0|>x\right)}\\
&\stw& \P\left( \left(\left(\barr{rr} \bfTh_0\\Y_{p,0}|\bfTh_0|\\|\bfTh_0|\earr\right)\,,\ldots,
\left(\barr{rr} \bfTh_h\\Y_{p,h}|\bfTh_h|\\|\bfTh_h|\earr\right) \right)\in\cdot \right)\,,\qquad \xto\,.
\eeao
The limit highlights the \regvar\ structure of $(\bfZ_t)$.\\[2mm]
The anti-clustering condition \eqref{cond:acac} is easily checked on $(\bfZ_t)$ by exploiting the corresponding properties of $(\bfX_t)$ and using a domination argument 
on the conditional expectation given $(Y_{p,t})$.
Exploiting the alternative mixing condition \eqref{eq:wdepms2}, 
we may proceed as in the proof of Theorem~\ref{th:joint} to conclude that, as $\nto$,
\beam\label{eq:convlog}\lefteqn{
k_n\log\big(\E\big[\ex^{i \,a_n^{-1}\,\bfu^\top\,\bfS_{r_n}}\,\1\big(a_n^{-1}\wt M_{r_n}\le \la^{-1/p}\,, a_n^{-1}M_{r_n}^{|\bfX|}\le x\big)\big]\big)}\\[2mm]
&\to&\E[\log\varphi_{\bfxi_\a}(\bfu)]\nonumber\\&&-
\dint_0^\infty\E\Big[\ex^{i\,y\, \sum_{t\in\bbz} \bfu^\top\bfQ_t}\1\Big(y\,\max_{t\in \bbz}\big((Y_{p,t}\,\la^{1/p})\vee (1/x)\big)|\bfQ_t|> 1 \Big)\Big]\,d(-y^{-\a})\nonumber\\[2mm]
%\eeao
%where
%\beao
%\lefteqn{\E[\log\varphi_{\bfxi_\a}(\bfu)]-
%\dint_0^\infty\E\Big[\ex^{i\,y\, \sum_{t\in\bbz} \bfu^\top\bfQ_t}\1\Big(y\,\max_{t\in \bbz}\big((Y_{p,t}/\la^{1/p}\vee (1/x)\big)|\bfQ_t|> 1 \Big)\Big]\,d(-y^{-\a})}\\
&=& \int_0^\infty 
\E\Big[\ex^{i\,y\, \sum_{t\in\bbz} \bfu^\top \bfQ_t}\,
\1\Big(y\,\max_{t\in \bbz}\,Y_{p,t}\,|\bfQ_t|\le \la^{-1/p}\,,
y\,\max_{t\in \bbz}\,|\bfQ_t|\le x \Big)
\nonumber\\&& \hspace{0.9cm}-1 - 
i\,y\,\sum_{t\in\bbz}\, \bfu^\top\bfQ_t\,\1_{(1,2)}(\a) \Big]\,d(-y^{-\a})\nonumber\\
&=& \int_0^\infty 
\E\Big[\ex^{i\,y\, \bfu^\top\sum_{t\in\bbz}\bfQ_t\ -y^p\,\la\,\sum_{t\in\bbz}|\bfQ_t|^p}\,
\1\Big(y\,\max_{t\in \bbz}|\bfQ_t|\le x \Big)\nonumber\\ &&\hspace{0.9cm} -1 - 
i\,y\,\bfu^\top \sum_{t\in\bbz} \,\bfQ_t\,\1_{(1,2)}(\a) \Big]\,d(-y^{-\a})\,.\nonumber
\eeam
In the last identity we again used the \ds\ of the iid Fr\'echet variables $Y_{p,t}$, $t\in \bbz$.
The desired result follows.
\end{proof}

\section{Ratio limit theorems for self-normalized quantities}\label{sec:ratiolmit}\setcounter{equation}{0}
\subsection{Ratios of sums and maxima}\label{subsec:ratio}
A direct con\seq\ of Theorem~\ref{th:joint} and the   \cmt\  is the following 
characterization of the limit distribution of the ratios
$(\bfS_n/M_n^{| \bfX|})$.

\bco\label{cor:ratiolimit}
Under the assumptions of Theorem \ref{th:joint}  we have that
\beao
\dfrac{\bfS_n}{M_n^{|\bfX|}}\std \dfrac{\bfxi_\a}{\eta_\a}=:\bfR_\a\,,\qquad \nto\,,
\eeao
where the random vector  $\bfR_\a$ has characteristic function, for every $\bfu\in\bbr^d$,
\beao
\varphi_{{\bfR}_\a}(\bfu)=\dfrac{ \E\big[\ex^{i\, \bfu^\top\sum_{t=-\infty}^\infty \wt  \bfQ_t}  \big]}{ \dint_0^\infty\E\Big[1 + 
i\,y\,\bfu^\top\,\sum_{t\in\bbz}\wt\bfQ_t\,\1_{(1,2)}(\a) -\ex^{iy \bfu^\top \sum_{t=-\infty}^\infty \wt\bfQ_t}\1(y\le1)  \Big]\,d(-y^{-\a})} \,.
\eeao 
\eco
%The random vector ${\bfR}_\a$ has the same moments properties than $\sum_{t=-\infty}^\infty \wt  \bfQ_t$. 
\bre
\label{rem:existence:sumQ}
In the expression for $\varphi_{\bfR_\a}$ we need to ensure that $\sum_{t\in\bbz} \wt\bfQ_t$ is well defined. If $1<\a<2$ and the
anti-clustering condition \eqref{cond:acac} holds we have $
\E\Big[\Big(\sum_{t\in\bbz}|\bfQ_t|\Big)^\a\Big] <\infty$; see \eqref{eq:xmasjan4a}.
By the definition of $(\wt\bfQ_t)$ (cf. \eqref{eq:qtilde}), the fact that 
$\max_{t\in \bbz}|\bfQ_t|\le 1$
 and Jensen's inequality we get 
\beao
 \E\Big[\sum_{t\in\bbz}  |\wt\bfQ_t|\Big]&=&\theta_{|\bfX|}^{-1}\E\Big[\big(\max_{t\in \bbz}|\bfQ_t|\big)^{\a-1}\sum_{t\in\bbz}  |\bfQ_t|\Big]
\le\theta_{|\bfX|}^{-1}\E\Big[\sum_{t\in\bbz}  |\bfQ_t|\Big]\\
&\le&\theta_{|\bfX|}^{-1} \Big(\E\Big[\Big(\sum_{t\in\bbz}  |\bfQ_t|\Big)^\a\Big]\Big)^{1/\a}<\infty\,.
\eeao
For  $0<\a<1$ 
similar calculations yield
\beao
 \E\Big[\sum_{t\in\bbz}  |\wt\bfQ_t|\Big]=\theta_{|\bfX|}^{-1}\E\Big[\dfrac{\sum_{t\in\bbz} | \bfTh_t|}{\big(\max_{t\in \bbz}|\bfTh_t|\big)^{1-\a}\sum_{t\in\bbz}  |\bfTh_t|^\a}\Big]
\le\theta_{|\bfX|}^{-1}<\infty\,.
\eeao
We can also derive
\beam\label{eq:ralpha}
\E[{\bfR}_\a]=\dfrac{\E\Big[\sum_{t\in\bbz} \wt \bfQ_t\Big]}{1-\a}\,,
\eeam
see Appendix \ref{app:domratiobis}. %{app:domratio}. 
The \rhs\ term is well
defined since $\sum_{t\in\bbz}  |\wt\bfQ_t|$ has finite expectation.
For iid centered
\rv s $(X_i)$ this is 
in agreement with the known fact  that
$\E[{R}_\a]=(1-\a)^{-1}=\lim_{\nto}
\E[S_n/M_n]$; see Bingham et al. \cite{bingham:goldie:teugels:1987}. The convergence of 
$(\E[\bfS_n/M_n^{|\bfX|}])$  for dependent $(\bfX_t)$ is an open question.
\ere

\begin{proof} By virtue of Theorem~\ref{th:joint} and the \cmt\
the limit relation $\bfS_n/\bfM_n^{|\bfX|}\std \bfR_\a$ follows. It remains 
to identify the \chf\ of~$\bfR_\a$.
\par
We re-write \eqref{eq:maxsum} as follows: for $x>0$ and $\bfu\in\bbr^d$,
\beao
\Psi(\bfu,x)&=&\E\big[\ex^{i\bfu^\top \bfxi_\a}\1( \eta_\a\le x)\big]\\
&=&\varphi_{\bfxi_\a}(\bfu)\exp\Big(-\theta_{|\bfX|}\dint_x^\infty\E\big[\ex^{i\,y\,\bfu^\top\sum_{t=-\infty}^\infty \wt \bfQ_t}  \big]\,d(-y^{-\a})\Big)\,.
\eeao
Denote the density of $\eta_\a$ by $f_{\eta_\a}$. We differentiate both 
sides of the latter identity \wrt~$x$:
\beao
\E\big[\ex^{i\,\bfu^\top \bfxi_\a}\;\big|\; \eta_\a= x\big]\,f_{\eta_\a}(x)
&=&\Psi(\bfu,x)\,\,\a\, x^{-\a-1}\,\theta_{|\bfX|}\,\underbrace{\E\big[\ex^{i\,x\,\bfu^\top\sum_{t=-\infty}^\infty \wt \bfQ_t}  \big]}_{ =: g(x\,\bfu)}\,.
%&=:&\Psi(\bfu,x)\,\,\a\, x^{-\a-1}\,\theta_{|\bfX|}\,g(x\bfu)\,.
\eeao
In particular, we obtain 
\beao
%\E\big[\ex^{i\bfu^\top \bfxi_\a/\eta_\a}\;\big|\; \eta_\a= x\big]\,f_{\eta_\a}(x)
\E\big[\ex^{i\bfu^\top \bfxi_\a/x}\;\big|\; \eta_\a= x\big]\,f_{\eta_\a}(x)
%&=&\Psi(\bfu/x,x)\,\,\a\, x^{-\a-1} \theta_{|\bfX|}\,\E\big[\ex^{i\, \bfu^\top\sum_{t=-\infty}^\infty \wt \bfQ_t}  \big]\\
&=&\Psi(\bfu/x,x)\,\,\a\, x^{-\a-1} \theta_{|\bfX|} \,\,g(\bfu)\,.
\eeao
Integration yields the following expression for $\varphi_{\bfR_\a}(\bfu)$:
\beao
&&\int_0^\infty\E\big[\ex^{i\bfu^\top \bfxi_\a/\eta_\a}\;\big|\; \eta_\a= x\big]\,f_{\eta_\a}(x) \,dx\\
&=&g(\bfu)\,\theta_{|\bfX|} \,\int_0^\infty\Psi(\bfu/x,x)\,d(-x^{-\a})\\
&=&g(\bfu)\,\theta_{|\bfX|} \,\int_0^\infty\varphi_{\bfxi_\a}(\bfu/x)\,\exp\Big(-\theta_{|\bfX|}\dint_x^\infty\E\big[\ex^{i\,(y/x)\,\bfu^\top\sum_{t=-\infty}^\infty \wt \bfQ_t}  \big]\,d(-y^{-\a})\Big)\,d(-x^{-\a})\\
&=&g(\bfu)\,\theta_{|\bfX|} \,\int_0^\infty\exp\Big(x^{-\a}\Big(\log \varphi_{\bfxi_\a}(\bfu)-\theta_{|\bfX|}\dint_1^\infty \E\big[\ex^{i\,y\,\bfu^\top\sum_{t=-\infty}^\infty \wt \bfQ_t}  \big]\,d(-y^{-\a})\Big)\Big)\,d(-x^{-\a})\,.
\eeao
Plugging in the \chf\ $\varphi_{\bfxi_\a}$ and 
changing variables, we obtain
\beao
&&g(\bfu)\,\theta_{|\bfX|} \,\int_0^\infty\exp\Big(z\theta_{|\bfX|} \dint_0^\infty\E\Big[\ex^{iy \bfu^\top \sum_{t=-\infty}^\infty \wt\bfQ_t}\1(y\le 1)  -1 - 
i\,y\,\bfu^\top\,\sum_{t\in\bbz}\wt\bfQ_t\,\1_{(1,2)}(\a) \Big]\\&&\hspace{9cm}d(-y^{-\a})\Big)\,dz\\
&&=\dfrac{ \E\big[\ex^{i\, \bfu^\top\sum_{t=-\infty}^\infty \wt  \bfQ_t}  \big]}{ \dint_0^\infty\E\Big[1+ 
i\,y\,\bfu^\top\,\sum_{t\in\bbz}\wt\bfQ_t\,\1_{(1,2)}(\a)- \ex^{iy \bfu^\top \sum_{t=-\infty}^\infty \wt\bfQ_t}\1(y\le1) \Big]\,d(-y^{-\a}) }\,.
\eeao
This is the desired formula for the \chf\ of $\bfR_\a$. 
\end{proof}
\bre\label{rem:regvary}
By the \cmt\ and a similar argument as in the proof above we have
\beao
\lim_{m\to\infty}\lim_{\nto}  \E\big[\ex^{i\bfu^T \bfS_n/M_n^{|\bfX|}}\mid M_n^{|\bfX|}>m\,a_n\big]
&=&\lim_{m\to\infty}\E\big[\ex^{i\bfu^T \bfR_\a}\mid \eta_\a>m\big]\\
&=&
\lim_{m\to\infty}\dfrac{g(\bfu)\int_m^\infty\Psi(\bfu/x,x)\,d(-x^{-\a})}{g({\bf0})\int_m^\infty\Psi({\bf0},x)d(-x^{-\a})}\\&=& 
g(\bfu)\lim_{m\to\infty}\dfrac{\int_0^{m^{-\a}}\Psi(\bfu\,y^{1/\a},y^{-1/\a})\,dy}{ \int_0^{m^{-\a}}\Psi({\bf0},y^{-1/\a})dy}\\
&=&
g(\bfu)\lim_{m\to\infty}\dfrac{\Psi(\bfu/ m,m)}{\Psi({\bf0},m)}\\
&=&g(\bfu)\,.
\eeao
Here we also used l'Hospital's rule and the explicit form of the 
hybrid \chf\ $\Psi$, satisfying $\Psi(\bfu/m,m)\to 1$ as $m\to \infty$.
Because $g(\bfu)= \E\big[\ex^{i\, \bfu^\top\sum_{t=-\infty}^\infty \wt  \bfQ_t}  \big]$, $\bfu\in \R^d$, this implies that as $m\to\infty$,
\beao
\lim_{\nto}\P\Big( \dfrac{\bfS_n}{M_n^{|\bfX|}}\in\cdot \Big| M_n^{|\bfX|}>m\,a_n\Big)
&=&\P\Big(\dfrac{\bfxi_\a}{\eta_\a}\in\cdot \Big|\eta_\a>m\Big)
\stw \P\Big(\sum_{t=-\infty}^\infty \wt  \bfQ_t\in \cdot\Big)\,.
\eeao
The case $\bfxi_\a=\bf0$ corresponds to $\sum_{t=-\infty}^\infty \wt  \bfQ_t=\bf0$ a.s. In what follows, we assume that this condition is not satisfied.
Then $|\bfxi_\a|$ and $\eta_\a$ are tail-equivalent.
Using the \regvar\ property of $\eta_\a$, we get the convergence of 
\beao
\P\big(m^{-1}(\bfxi_\a,\eta_\a)\in \{(\bfu,x)\in\R^d\times (0,\infty): (\bfu/x,x)\in 
A\times (y,\infty)\}\big)/\P(\eta_\a> m )\,,\qquad m\to\infty\,,
\eeao 
for every continuity set $A$ of the 
distribution of $\sum_{t=-\infty}^\infty \wt  \bfQ_t$ and $y>0$. The sets $A\times (y,\infty)$ generate the vague \con\ in $\R^d\times(0,\infty)$ of these \ms s, implying  
that  $(\bfxi_\a,\eta_\a)$ is an $\bbr^{d}\times\R_+$-valued 
\regvary\ vector with index $\a>0$. 
The \regvar\ properties of   $(\bfxi_\a,\eta_\a)$ then follows from the vague convergence on $\R^d\times(0,\infty)$ extended to $(\R^{d}\times\R_+)\setminus\{{\bf0}\}$, using the \regvar\ properties  of the $\a$-stable distribution of $\bfxi_\a$.
\ere
\begin{comment}
{\red We illustrate the \ds\ of $R_\a$ in the case of a \regvary\ stationary
AR(1) process.}
\begin{figure}[thbp]
\centerline{
\epsfig{figure=ratios05iid.eps,width=6cm,height=6cm,angle=-90}\epsfig{figure=ratios15iid.eps,width=6cm,height=6cm,angle=-90}}
\centerline{
\epsfig{figure=ratios05ar09.eps,width=6cm,height=6cm,angle=-90}\epsfig{figure=ratios15ar09.eps,width=6cm,height=6cm,angle=-90}}
\centerline{
\epsfig{figure=ratios05ar-09.eps,width=6cm,height=6cm,angle=-90}\epsfig{figure=ratios15ar-09.eps,width=6cm,height=6cm,angle=-90}}
\caption{  Boxplots based on  $10^3$ Monte-Carlo simulations of 
$n=10^5$ observations of an AR(1) model $X_t= \varphi\,X_{t-1}+Z_t$ for 
iid Pareto$(\a)$ \rv s $(Z_t)$. {\bf Left:} The \ds\ of the
ratios $S_n/M_n^{|X|}$ for $\alpha=.5$. {\bf Right:}
The \ds\ of the ratios  $(S_n-n\,\E[X])/M_n^{|X-\E[X]|}$ for $\a=1.5$. 
{\bf From top to bottom:} $\varphi=0, -0.9,0.9$. The {\blue blue} and {\red red} lines
indicate the sample mean of the ratio and $\E[R_\a]$, respectively. The simulations are in good agreement with the theoretical value 
of~$\E[R_\a]=(1-\varphi)/(1-\a)$.
}\label{fig:ratios}
\end{figure}
\end{comment}
\subsection{Studentized sums}\label{subsec:student}
Corollary~\ref{cor:ratiolimit} deals with a special case of self-normalization 
of the (possibly centered) sum processes $(\bfS_n)$ constructed from an
$\bbr^d$-valued \regvary\ \seq\ $(\bfX_t)$ with index $\a\in (0,2)\backslash\{1\}$. Indeed,
by virtue of the joint \con\ of the partial sums and maxima 
in Theorem~\ref{th:joint}
with the same
normalization $(a_n)$ the latter \seq\  can be replaced by $(M_n^{|\bfX|})$, 
leading to a ratio limit theorem for sums and maxima. An advantage of 
this approach is that the {\em unknown} \seq\ $(a_n)$ is replaced by
the observed maxima $(M_n^{|\bfX|})$. The price one has to pay for this is that
the limit \ds\ of the ratio is an unfamiliar \ds\ which does not belong to the $\a$-stable class and, for its evaluation,  one has to employ Monte-Carlo simulations or numerical techniques.
\par
The classical example of self-normalization is studentization: 
the standard deviation in the normalization of the classical \clt\
is replaced by an empirical estimator. Our goal is to derive limit theory
for the partial sums $(\bfS_n)$ with the $\ell^p$-norm-type moduli
$\gamma_{n,p}$ for $p>0$ defined in \eqref{eq:moduli}. 
We assume $\a<2$ and $p>\a $. Then $(|\bfX_t|^p)$ is a \regvary\ \seq\
with index $\a/p$. In particular, we can always choose $p=2$, corresponding
to the classical studentization.
\par
Following the ideas of the proof of Corollary~\ref{cor:ratiolimit}, 
we obtain a limit for the ratios $(\gamma_{n,p}/M_n^{|\bfX|})$.
\bco
Under the assumptions of Theorem~  %\ref{th:joint}  
 \ref{pr:prophybridch} we have 
\beao
\dfrac{\gamma_{n,p}}{M_n^{|\bfX|}}\std \dfrac{\zeta_{\a,p}}{\eta_\a}\,,\qquad \nto.
\eeao
For $p>\a$, $(\zeta_{\a,p}/\eta_\a)^p$ has Laplace transform
\beao
\E[\ex^{-\la\, (\zeta_{\a,p}/\eta_\a)^p}]&=&\dfrac{ \E\big[\ex^{-\la\,\sum_{t=-\infty}^\infty |\wt  \bfQ_t|^p}  \big]}{ \dint_0^\infty\E\big[1-\ex^{-y^p\,\la\,  \sum_{t=-\infty}^\infty |\wt\bfQ_t|^p}\1(y\le1)  \big]\,d(-y^{-\a})} \,,\qquad \la>0\,.
\eeao 
\eco
The proof follows the lines of the one for Corollary \ref{cor:ratiolimit}; it 
is omitted. Similarly to Remark~\ref{rem:regvary} we can also derive that 
$(\zeta_{\a,p},\eta_\a)$ is an $\R_+^2$-value \regvary\ vector with index $\a$, satisfying 
\beao
 \P\Big(\dfrac{\zeta_{\a,p}}{\eta_\a}\in \cdot \Big| \eta_\a>m\Big)\stw \P\Big(\sum_{t=-\infty}^\infty |\wt  \bfQ_t|^p\in\cdot\Big),\qquad m\to\infty\,.
\eeao
The studentized sums  also converge
\beam\label{eq:mod}
\dfrac{\bfS_n}{\gamma_{n,p}}\std \dfrac{\bfxi_\a}{\zeta_{\a,p}}=:{\bfR}_{\a,p}\,,\qquad \nto\,.
\eeam
However, it is even more difficult to describe 
the limit distribution. We can still derive the first moment of the ratio 
${\bfR}_{\a,p}$ from the joint \chf -- Laplace transform 
$\Phi_{\bfxi_\a,\zeta_{\a,p}^p}(\bfu,\lambda):=
\E\big[\ex^{i\,\bfu^\top \,\bfxi_\a-\lambda  \,\zeta_{\a,p}^p}\big]$, 
$(\bfu,\lambda)\in\bbr^d\times\R_+$, given in Theorem~\ref{pr:prophybridch}.
\bpr\label{pr:ratiostud}
Under the assumptions of Theorem~\ref{pr:prophybridch}, \eqref{eq:mod} holds
and
\beam\label{eq:momentgoodexp}
\E[\bfR_{\a,p}]=\dfrac{\Gamma((1-\a)/p)}{\Gamma(1/p)\Gamma(1-\a/p)} \E\Big[\frac{\|\bfQ\|_p^\a}{\E[\|\bfQ\|_p^\a]}\dfrac{\sum_{t=-\infty}^\infty\bfQ_t}{\|\bfQ\|_p}\Big]\,.
\eeam
\epr
\begin{proof}
We start by verifying  the finiteness of $\E[\|\bfQ\|_p^{\a-1}\|\bfQ\|_1]$, $ \a<p$, ensuring 
the integrability properties needed throughout this proof. 
For $\a\ge 1$  we use the monotonicity of the $\ell^p$-norms 
$\|\bfQ\|_p\le \|\bfQ\|_1$ for $p\ge 1$. 
We obtain $\E[\|\bfQ\|_p^{\a-1}\|\bfQ\|_1]\le \E[\|\bfQ\|_1^\a]$ which is finite; see Remark \ref{rem:expnormfinite}. 
For $\a<1$ we have $\|\bfQ\|_p\ge \|\bfQ\|_\infty$ hence
\beao
\|\bfQ\|_p^{\a-1}\|\bfQ\|_1\le \|\bfQ\|_\infty^{\a-1}\|\bfQ\|_1\le \|\bfQ\|_\a^\a= 1\,,
\eeao
by definition of the spectral cluster process $(\bfQ_t)$.
\par
The following formula  is immediate from the definition of the $\Gamma$-function:
\beam\label{eq:gammaformula}
\dfrac1{x^{1/p}}=\dfrac  1{\Gamma(1/p)}\int_0^\infty\la^{1/p-1}\ex^{-\la x}d\lambda=\dfrac  p{\Gamma(1/p)}\int_0^\infty\ex^{-\la^p x}d\lambda\,,\qquad x>0\,.
\eeam
We need the following lemma whose proof is postponed to Appendix \ref{app:domratio}.
\ble\label{lem:domratio}
Under the assumptions of Theorem~\ref{pr:prophybridch} we have $\E\big[|\bfxi_{\a}|\ex^{-\lambda^p\zeta_{\a,p}^p}\big]<\infty$.
\ele
Then an application of Fubini's theorem and \eqref{eq:gammaformula} yields
\beao
\E[\bfR_{\a,p}]&=&\E\Big[\dfrac{\bfxi_{\a}}{(\zeta_{\a,p}^p)^{1/p}}\Big] 
=\dfrac  p{\Gamma(1/p)}\int_0^\infty \E\big[\bfxi_{\a}\ex^{-\lambda^p\zeta_{\a,p}^p}\big]d\lambda\,.
\eeao
By another application of Lemma \ref{lem:domratio} the integrand coincides with 
\beao
\dfrac1i\dfrac{\partial \Phi_{\bfxi_\a,\zeta_{\a,p}^p}({\bf0},\lambda^p)}{\partial \bfu}&=&\int_0^\infty\E\Big[y\,\sum_{t=-\infty}^\infty\bfQ_t\,
\big(\ex^{-(\la y)^p\|\bfQ\|_p^p}-\1_{(1,2)}(\a)\big)
\Big]\,d(-y^{-\a})\;\E\big[\ex^{-\lambda^p\zeta_{\a,p}^p}\big]\,.
\eeao
Exploiting  \eqref{eq:gammaformula} for $0<\alpha<1$, the formula 
\beao
\int_0^\infty z\, \big(\ex^{-z^p\, x}-1\big)\, d (-z^{-\a})=
(\alpha/p)\Gamma((1-\alpha)/p)\,x^{(\alpha-1)/p}\,,
\eeao
for $\a>1$,
and changing the variable to $z=\la\, y\, \|\bfQ\|_p$, we obtain
\beao\lefteqn{
\int_0^\infty\E\Big[y\sum_{t=-\infty}^\infty\bfQ_t\ex^{-(\la y)^p\|\bfQ\|_p^p}
\big(\ex^{-(\la y)^p\|\bfQ\|_p^p}-\1_{(1,2)}(\a)\big)
\Big]d(-y^{-\a})}\\&=&\frac\a p\Gamma((1-\a)/p)\la^{\a-1}\E\Big[\|\bfQ\|_p^\a\dfrac{\sum_{t=-\infty}^\infty\bfQ_t}{\|\bfQ\|_p}\Big]\,.
\eeao
Combining these results, one obtains
\beao
\E[\bfR_{\a,p}]&=&\dfrac  {\Gamma((1-\a)/p)}{\Gamma(1/p)}\E\Big[\|\bfQ\|_p^\a\dfrac{\sum_{t=-\infty}^\infty\bfQ_t}{\|\bfQ\|_p}\Big]\int_0^\infty \a\la^{\a-1}\E\big[\ex^{-\lambda^p\zeta_{\a,p}^p}\big]d\la\,.
\eeao
We recall the Laplace transform of $\zeta_{\a,p}^p$ derived in Remark \ref{rem:stableskew}, i.e.,
\beao
\E\big[\ex^{-\lambda^p\zeta_{\a,p}^p}\big]=\exp\big(-\Gamma(1-\a/p)\,\E[\|\bfQ\|_p^\a]\,\lambda^\a\big)\,,\qquad\lambda>0\,.
\eeao
Plugging in the previous expression, we arrive at
\beao
\E[\bfR_{\a,p}]&=& \dfrac  {\Gamma((1-\a)/p)}{\Gamma(1/p)}\E\Big[\|\bfQ\|_p^\a\dfrac{\sum_{t=-\infty}^\infty\bfQ_t}{\|\bfQ\|_p}\Big] \int_0^\infty\ex^{-\Gamma(1-\a/p)\E[\|\bfQ\|_p^\a]\lambda^\a}d(\la^\a)\\
&=&\dfrac  {\Gamma((1-\a)/p)}{\Gamma(1/p)}\E\Big[\|\bfQ\|_p^\a\dfrac{\sum_{t=-\infty}^\infty\bfQ_t}{\|\bfQ\|_p}\Big]\dfrac1{\Gamma(1-\a/p)\E[\|\bfQ\|_p^\a]}\,.
\eeao
The desired result follows. 
\end{proof}

\subsubsection*{The case $p=2$}
It will be convenient to introduce the spectral cluster process $\wh\bfQ$ of the \regvary\ \seq\ $(\bfX_t)$ 
by a change of measure
\beao
\P(\wh \bfQ\in \cdot)=\E\Big[\dfrac{\|\bfQ\|_2^\a}{\E[\|\bfQ\|_2^\a]}\1\Big(\dfrac{\bfQ}{\|\bfQ\|_2}\in \cdot\Big)\Big]\,.
\eeao
Then \eqref{eq:momentgoodexp} turns into
\beao
\E[\bfR_{\a,2}]=\dfrac{\Gamma((1-\a)/2)}{\Gamma(1/2)\Gamma(1-\a/2)} \E\Big[ \sum_{t=-\infty}^\infty\wh\bfQ_t \Big]\,.
\eeao
\subsection{Greenwood statistics}\label{subsec:greenwood}
We consider a univariate positive \regvary\ stationary \ts\ $(X_t)$ 
with index $\a \in (0,1)$.  Its spectral tail process is denoted 
by $(\Theta_t)$ and the corresponding spectral cluster process by
$Q=(Q_t)=(\Theta_t/\|\Theta\|_\a)$. The {\em Greenwood statistic} (see Greenwood \cite{greenwood:1946}) is the ratio statistic
\beao
\dfrac{\gamma_{n,2}^2}{S_n^2}= \dfrac{X_1^2+\cdots + X_n^2}
{(X_1+\cdots +X_n)^2}\,, \qquad n\ge 1\,.
\eeao
Under the conditions of Theorem~\ref{pr:prophybridch} 
a continuous mapping argument yields the more general result for $\a<p$  
\beao
T_{n,p}:=\dfrac{X_1^p+\cdots + X_n^p}
{(X_1+\cdots +X_n)^p}
\std \dfrac{\zeta_{\a,p}^p}{\xi_\a^p}\,,
\eeao
where $\zeta_{\a,p}^p$ is totally skewed to the right $\a/p$-stable and 
 $\xi_\a$ is $\a$-stable. %In contrast to the studentized statistics we observe 
%that the \lhs\ is bouded by 1 and dominated \con\ yields  
%\beao
%\E[T_{n,p}] \to \E \Big[ \dfrac{\zeta_{\a,p}^p}{\xi_\a^{p}}\Big]
%\eeao
\bco\label{cor:greenwood}
Under the conditions of  Theorem~\ref{pr:prophybridch}  we have 
for $\a<p\wedge 1$,
\beao
\E[T_{n,p}] \to \E \Big[ \dfrac{\zeta_{\a,p}^p}{\xi_\a^{p}}\Big]
%\E \Big[ \dfrac{\zeta_{\a,p}}{\xi_\a^{p}}\Big]
&=&
\frac{\Gamma(p-\alpha)}
%\Gamma(\frac{q-p}{\alpha}+1)}
{\Gamma(p)\,\Gamma(1-\alpha)}
\, \E\Big[\frac{
\|Q\|_1^\alpha}{\E[\|Q\|_1^\alpha] } 
\frac{\|Q\|_p^p}{\|Q\|_1^p}
\Big]\,.
%\E [\|Q\|_1^\alpha]^{\frac{q-p}{\alpha}}.
\eeao 
In particular,
\beao
\E \Big[ \dfrac{\zeta_{\a,2}^2}{\xi_\a^{2}}\Big]
&=&
(1-\alpha)\, \E\Big[\frac{
\|Q\|_1^\alpha}{\E[\|Q\|_1^\alpha]} 
\frac{\|Q\|_2^2}{\|Q\|_1^2}
\Big].
\eeao 
\eco
In the case of \asy\ independence when $\Theta_t=0$ for $t\ne 0$ the 
\rhs\ turns into $\Gamma(p-\alpha)/(\Gamma(p)\,\Gamma(1-\alpha))$, in agreement with Albrecher et al. \cite{albrecher:teugels:2006,albrecher:flores:2022} in the iid case. 
\begin{proof}
Considering $X_t'=X_t^p$ we obtain a \regvary\ \seq\ of index $\a/p$ with spectral cluster process $(Q_t^p)$. The limit can be interpreted as the ratio $R_{\a/p,1/p}'=\xi_{\a/p}'/\zeta_{\a/p,1/p}'$ with joint \chf --Laplace transform
of  $(\xi_{\a/p}',\zeta_{\a/p,1/p}')$ given by
\beao
\E\big[\ex^{i \,u\, \xi_{\a/p}'-\la\, (\zeta_{\a/p,1/p}')^{1/p}}\big]
&=&\exp\Big(\dint_0^\infty\E\Big[\ex^{i\,y\, u\|Q\|_p^p-y^{1/p}\la\,\|Q\|_1}
-1\Big]d(-y^{-\a/p})\Big)\,,\qquad u\in\R,\,\lambda>0\,.
%&=&\varphi_{\bfxi_\a}(\bfu)\,
%\exp\Big(-\dint_0^\infty\E\Big[\ex^{i\,y\,\bfu^\top \sum_{t=-\infty}^\infty Q_t}\, \1\Big(y\,\max_{t\in \Z}|Q_t|>x \Big)\Big]\,d(-y^{-\a})\Big)\nonumber\\
%&=&\varphi_{\bfxi_\a}(\bfu)\,\Phi_\a^{\theta_{|\bfX|}}(x)\, \exp\Big(-\theta_{|\bfX|}\dint_x^\infty\E\big[\ex^{i\,y\,\bfu^\top\sum_{t=-\infty}^\infty \wt \bfQ_t}-1\big]\,d(-y^{-\a})\Big)\nonumber\\
%&=&\varphi_{\bfxi_\a}(\bfu)\,\Phi_\a^{\theta_{|\bfX|}}(x)\,
%\exp\big(-\E\big[\ex^{i\,\bfu^\top Y\,\sum_{t=-\infty}^\infty \wt \bfQ_t}-1\big]\big)\,,\nonumber\\
\eeao
The expression of the first moment follows from an application of Proposition \ref{pr:ratiostud}.
\end{proof}
\bexam
\rm Consider the scaled sample kurtosis for 
the stationary sequence $(X_t)_{t \in \mathbb{Z}}$:
\beao
\frac{\sum_{i=1}^n |X_t|^4}{(\sum_{i=1}^n |X_t|^2 )^2} = \frac{\|X\|_4^4}{\|X\|^4_2} \le 1\,.
\eeao
Assume that $(X_t)$ is \regvary\ with index $0<\a<2$. Then $(|X_t|^2)$ is also \regvary\ with index $0<\a/2<1$ and if it satisfies the assumptions of Corollary \ref{cor:greenwood} we obtain, keeping the same notation,
\beao
\frac{\|X\|_4^4}{\|X\|^4_2} \stackrel{d}{\to} \frac{\zeta_{\alpha/2,2}^2}{ \xi_{\alpha/2}^2} \quad \text{and} 
\quad \E\Big[\frac{\|X\|_4^4}{\|X\|^4_2}\Big]\to \E\Big[ \frac{\zeta_{\alpha/2,2}^2}{ \xi_{\alpha/2}^2} \Big]=(1-\a/2)\E\Big[\frac{\|Q\|_2^\a}{\E[\|Q\|_2^\a]}\frac{\|Q\|_4^4}{\|Q\|_2^4}\Big]. 
\eeao
\eexam

\subsection{Ratios of norms}
\label{subsec:ration:norms}
Consider a \regvary\ stationary process $(\bfX_t)$ 
with index $\a>0$. For positive $q>0$ we consider the 
norm-type modulus of the sample $\bfX_1,\ldots,\bfX_n$ given by
$\|\bfX\|_q=(\sum_{t=1}^n|\bfX_t|^q)^{1/q}$. Here we suppress the dependence
on $n$ in the notation.
\par
For $q,r>0$ we  are interested in the limit behavior of the ratios
${\|\bfX\|_q}/{\|\bfX\|_r}$.
We rephrase this ratio in terms of the \regvary\ stationary \seq\ 
$(Z_t) = (|\bfX_t|^q)$ with index $\a/q$. By a  continuous mapping argument  
the tail process of $(Z_t)$ is given by $(|\bfTh_t|^q)$. Adapting the notation to the
 $Z$-\seq , we obtain
\beao
{\|\bfX\|_q}/{\|\bfX\|_r}=\big({\|\bfZ\|_1}/{\|\bfZ\|_{r/q}}\big)^{1/q}\,.
\eeao
If $(Z_t)$ satisfies the conditions of Theorem~\ref{pr:prophybridch}, 
in particular $\a<q\wedge r$, then
the continuous mapping theorem yields 
\beao
{\|\bfX\|_q}/{\|\bfX\|_r}\std R_{\a/q,r/q}^{1/q}\eqd {\xi_{\a/q}^{1/q}}/{\zeta_{\a/q,r/q}^{1/q}}\,,\qquad \nto\,.
\eeao
Here $(\xi_{\a/q},\zeta_{\a/q,r/q}^{ r/q})$ have the joint \chf\--Laplace transform
\beao\lefteqn{
\E\big[\ex^{i\,x\,\xi_{\a/q}-\la\,\zeta_{\a/q,r/q}^{r/q}}\big]}\\
&=&\exp\Big(\dint_0^\infty\E\Big[\ex^{i\,y\,x \,\sum_{t=-\infty}^\infty  Q_t -\la\,y^{r/q}\, \sum_{t=-\infty}^\infty  Q_t^{r/q}} -1\Big]d(-y^{-\a/q})\Big)\,,\qquad (x,\la)\in \bbr\times \bbr_+\,,
\eeao
and $(Q_t)=(|\bfTh_t|^q/\|\bfTh\|_q^\a)$ %\textcolor{red}{$(Q_t)=(|\bfTh_t|^q/\|\bfTh\|_\a^q)$ ?} 
is the spectral cluster process of $(Z_t)$. We also observe that for $\a< r\le q $ we have 
$\|\bfX\|_q/\|\bfX\|_r\le 1$ a.s., hence uniform integrability yields the \con\
of the moments
\beao
\E\big[{\|\bfX\|_q}/{\|\bfX\|_r}\big]\to 
\E\big[{\xi_{\a/q}^{1/q}}/{ \zeta_{\a/q,r/q}^{1/q}}\big]\,,\qquad \nto\,.
\eeao
It is desirable to derive a more explicit expression for the limit.

\section{Examples}\label{sec:examples}\setcounter{equation}{0}
\subsection{Sufficient conditions via coupling}
The anti-clustering condition \eqref{cond:acac} and the mixing condition \eqref{eq:wdepms2} can be checked by using a coupled version of $(\bfX_t)_{t\in \Z}$.
\bpr\label{prop:coupledmixed}  
We assume the following conditions:
\begin{enumerate}
\item[\rm (1)]
There exists a coupled version $(\bfX_t^\ast)$ of  
the \regvary\ stationary process $(\bfX_t)$ with index $\a\in(0,1)\cup (1,2)$: $(\bfX_t^\ast)$ is distributed as $(\bfX_t)$ and $(\bfX_t^\ast)_{t\ge 1}$ is independent of $(\bfX_t)_{t\le 0}$.
\item[\rm (2)]
For some integer \seq s  $(r_n)$ and $(\ell_n)$ such that  $r_n=o((a_n^{2}/n)\wedge n)$, $\ell_n=o(r_n)$, and for $q<\a\wedge1$, $p>\a$, we have
\beam\label{eq:couplingq}
k_n\,a_n^{-q}\sum_{t=\ell_n}^{r_n}\big(\E\big[|\bfX_t-\bfX_t^\ast|^q\big]\big)^{(1/p)\vee 1} \to 0\,,\qquad \nto\,.
\eeam
\item[\rm (3)]
For the same  \seq\ $(r_n)$ and the same $q<\a\wedge 1$ we have
\beam\label{eq:coupledcondition4}
\lim_{\kto}\limsup_{\nto} n \sum_{t=k}^{r_n}\E\big[\big(\big|a_n^{-1} \big(\bfX_t - \bfX_t^*\big)\big|^{q}\wedge 1\big) \big(\big|a_n^{-1}\bfX_0\big|^{q }\wedge 1\big)\big] = 0 \,. 
\eeam
\end{enumerate}
Then the anti-clustering condition \eqref{cond:acac} and the mixing condition  \eqref{eq:wdepms2} with $p$ as in \eqref{eq:couplingq}
hold for every $(\bfu,x,\la)\in \bbr^d\times \bbr_+^2$ and Theorem \ref{pr:prophybridch} applies.
\epr
\begin{proof}
We start by  checking \eqref{cond:acac}. Because $r_n=o((a_n^{2}/n)\wedge n)$ it is enough to show that \eqref{cond:centeracac} holds. Using the basic inequality 
\beao
|\cov(V,W)|\le \E[|W-W^\ast||V|]\,,
\eeao
where $W^\ast$ is a copy of $W$ independent of $V$ for $|V|\vee|W|\le 1$ a.s., we easily obtain the sufficient condition
\beao
\lim_{k\to\infty}\limsup_{\nto}n \,\sum_{j=k}^{r_n}\E\big[\big(\big|a_n^{-1}\big(\bfX_j - \bfX_j^*\big) \big|\wedge  1\big)(|  a_n^{-1}\bfX_{0}|\wedge
1)\big]=0 \,.
\eeao
The desired result follows by observing that $|x|\wedge 1\le |x|^{q'}\wedge 1$ for $0<q'\le 1$.
\par
Next we verify \eqref{eq:wdepms2}. We start by showing that
\eqref{eq:couplingq} implies 
\beam\label{eq:couplingq2}
k_n\,a_n^{-q}\sum_{t=\ell_n}^{r_n}\E\big[|\bfX_t-\bfX_t^\ast|^q+\big||\bfX_t|^p-|\bfX_t^\ast|^p\big|^{q/p}\big] \to 0\,,\qquad \nto\,.
\eeam
For $\a<p\le 1$ this follows by the triangular inequality
\beao
\big||\bfX_t|^p-|\bfX_t^\ast|^p\big|^{q/p}\le  |\bfX_t-\bfX_t^\ast|^{q}\,,\qquad t\ge1\,.
\eeao
For $p\ge 1$ we apply the mean value theorem:  there exists some $0<\xi<1$ such that 
\beao
\big||\bfX_t|^p-|\bfX_t^\ast|^p\big|^{q/p}=  \big(p||\bfX_t| + \xi(|\bfX_t^\ast| - |\bfX_t|)|^{(p-1)} ||\bfX_t|-|\bfX_t^\ast||\big)^{q/p}\,.
\eeao
Then the H\"older inequality yields
\beao
\E\big[\big||\bfX_t|^p-|\bfX_t^\ast|^p\big|^{q/p}\big]&\le& c_1\,\big(\E[||\bfX_t| + \xi(|\bfX_t^\ast| - |\bfX_t|)|^{q}\big]\big)^{(p-1)/p} \big(\E\big[|\bfX_t-\bfX_t^\ast|^{q}\big]\big)^{1/p}\\
&\le& c_2 \big(\E\big[|\bfX_t-\bfX_t^\ast|^{q}\big]\big)^{1/p}\,,\qquad t\ge 1\,,
\eeao
where $c_1$, $c_2>0$ do not depend on $t$. This proves \eqref{eq:couplingq2} for $p\ge 1$ as well.
\par
Next we introduce the intermediate \seq\ $(\ell_n)$ in the mixing condition \eqref{eq:wdepms2} by an application of an asymptotic negligibility argument on functionals on small blocks $(\bfX_{r_n(j-1)+t})_{1\le t\le \ell_n}$, $1\le j\le k_n$. 
Because the anti-clustering condition \eqref{eq:coupledcondition4} is also satisfied for the intermediate sequence $(\ell_n)$ we deduce similarly as in the proof of Theorem \ref{pr:prophybridch}, see equation \eqref{eq:convlog}, that the \seq s
\beao
\Big(\dfrac{n}{\ell_n}\log\big(\E\big[\exp\big(i\,a_n^{-1}\bfu^\top\bfS_{\ell_n}-a_n^{-p}\la \gamma_{\ell_n,p}^p\big)\,\1\big(a_n^{-1}M_{\ell_n}^{|\bfX|}\le x\big)\big]\big)\Big)
\eeao
converge for every $(\bfu,x,\la)$. Remembering that $\ell_n/r_n\to 0$ as $\nto$ we achieve that 
 \beao
\lefteqn{\lim_{\nto}\log\Big(\E\big[\exp\big(i\,a_n^{-1}\bfu^\top\bfS_{\ell_n}-a_n^{-p}\la \gamma_{\ell_n,p}^p\big)\,\1\big(a_n^{-1}M_{\ell_n}^{|\bfX|}\le x\big)\big]^{k_n}\Big)}\\
&=&\lim_{\nto}\dfrac{\ell_n}{r_n}\dfrac{n}{\ell_n}\log\big(\E\big[\exp\big(i\,a_n^{-1}\bfu^\top\bfS_{\ell_n}-a_n^{-p}\la \gamma_{\ell_n,p}^p\big)\,\1\big(a_n^{-1}M_{\ell_n}^{|\bfX|}\le x\big)\big]\big)=0\,.
\eeao
We immediately deduce the asymptotic neglibility of $k_n$ independent copies of the functionals on small blocks $(a_n^{-1}\bfS_{\ell_n},a_n^{-p} \gamma_{\ell_n,p}^p, a_n^{-1}M_{\ell_n}^{|\bfX|})$. Arguments similar to the ones developed after equation \eqref{eq:bigb} yield that
\beao
\lefteqn{\E\Big[\prod_{j=1}^{k_n}\exp\big(i\,a_n^{-1}\bfu^\top\bfS_{jr_n-\ell_n,jr_n}-a_n^{-p}\la \gamma_{jr_n-\ell_n,jr_n,p}^p\big)
\1\big(a_n^{-1}M_{jr_n-\ell_n,jr_n}^{|\bfX|}\le x\big)\Big]}\\
&-&
\Big(\E\big[\exp\big(i\,a_n^{-1}\bfu^\top\bfS_{\ell_n}-a_n^{-p}\la \gamma_{\ell_n,p}^p\big)\,\1\big(a_n^{-1}M_{\ell_n}^{|\bfX|}\le x\big)\big]\Big)^{k_n}=o(1)\,.
\eeao
We conclude the asymptotic negligibility of the functionals on small blocks $(\bfX_{r_n(j-1)+t})_{1\le t\le \ell_n}$, $1\le j\le k_n$ under the modified condition \eqref{eq:couplingq} such that $(\ell_n)$ is replaced by $(r_n-\ell_n)$.

Then the mixing condition \eqref{eq:wdepms2} turns into one based on large blocks only
\beam\label{eq:bigb}
\lefteqn{\E\Big[\prod_{j=1}^{k_n}\exp\big(i\,a_n^{-1}\bfu^\top\bfS_{j(r_n-1)+\ell_n,jr_n}-a_n^{-p}\la \gamma_{j(r_n-1)+\ell_n,jr_n,p}^p\big)
\1\big(a_n^{-1}M_{j(r_n-1)+\ell_n,jr_n}^{|\bfX|}\le x\big)\Big]}\\
&-&
\Big(\E\big[\exp\big(i\,a_n^{-1}\bfu^\top\bfS_{\ell_n,r_n}-a_n^{-p}\la \gamma_{\ell_n,r_n,p}^p\big)\,\1\big(a_n^{-1}M_{\ell_n,r_n}^{|\bfX|}\le x\big)\big]\Big)^{k_n}=o(1)\,.\nonumber
\eeam
We  use a telescoping sum argument over $1\le j\le k_n$ on the difference 
and show that it converges to zero as $\nto$.
Using the  properties of the coupled version,   this difference can be expressed as  a sum of $k_n$ summands plus negligible terms which we ignore. Up to a constant multiplier the absolute value of a typical summand
is bounded by
\beam
\E\big[\big|\ex^{i\,a_n^{-1}\,\bfu^\top\bfS_{\ell_n,r_n}-a_n^{-p}\,\la \,\gamma_{\ell_n,r_n,p}^p}\1\big(a_n^{-1}M_{\ell_n,r_n}^{|\bfX|}\le x\big)-\ex^{i\,a_n^{-1}\,\bfu^\top\bfS_{\ell_n,r_n}^\ast-a_n^{-p}\,\la\, \gamma_{\ell_n,r_n,p}^{\ast\,p}}\1\big(a_n^{-1}M_{\ell_n,r_n}^{|\bfX^\ast|}\le x\big)\big|\big]\,,\nonumber\\\label{eq:3x9}
\eeam
where $\bfS_{\ell_n,r_n}^\ast=\sum_{t=\ell_n}^{r_n}\bfX_t^\ast\,,$ $\gamma_{\ell_n,r_n,p}^{\ast\,p}=
\sum_{t=\ell_n}^{r_n}|\bfX_t^\ast|^p$ and $M_{\ell_n,r_n}^{|\bfX^\ast|}=\max_{\ell_n\le t\le r_n}|\bfX_t^\ast|$. 
We observe that \eqref{eq:3x9} is bounded by
\beao
\lefteqn{\E\big[\big|\ex^{i\,a_n^{-1}\bfu^\top\bfS_{\ell_n,r_n}} -
\ex^{i\,a_n^{-1}\bfu^\top\bfS_{\ell_n,r_n}^\ast}\big|\big]
+\E\big[\big|\ex^{-a_n^{-p}\,\la \,\gamma_{\ell_n,r_n,p}^p}-\ex^{-a_n^{-p}\,\la \,\gamma_{\ell_n,r_n,p}^{\ast \, p}}\big|\big]}\\
&&+\E\big[\big|\1\big(a_n^{-1}M_{\ell_n,r_n}^{|\bfX^|}\le x\big)-
\1\big(a_n^{-1}M_{\ell_n,r_n}^{|\bfX^\ast|}\le x\big)\big|\big]\\
&\le &
\E\big[\big|\big(a_n^{-1}\bfu^\top (\bfS_{\ell_n,r_n}-\bfS_{\ell_n,r_n}^\ast)\big)\wedge 1\big|^q\big]+
\E\big[\big|\big(a_n^{-p}\,\la \,(\gamma_{\ell_n,r_n,p}^p-\gamma_{\ell_n,r_n,p}^{\ast\, p})\big)\wedge 1\big|^{q'}\big]\\
&&+\Big[\P\big(a_n^{-1}M_{\ell_n,r_n}^{|\bfX|}> x\,,a_n^{-1}M_{\ell_n,r_n}^{|\bfX^\ast|}\le x\big)+
\P\big(a_n^{-1}M_{\ell_n,r_n}^{|\bfX^\ast|}> x \,,a_n^{-1}M_{\ell_n,r_n}^{|\bfX|}\le x\big)\Big]\\
&=:&I_1+I_2+I_3\,.
\eeao
Here we choose $q'=q/p <1$. Then
we have for some constant $c=c(\bfu,\la)>0$,
\beao
k_n(I_1+I_2)&\le & c\,k_n\,a_n^{-q}\,\sum_{t=\ell_n}^{r_n}\E\big[|\bfX_t-\bfX_t^\ast|^q\big]
+\big||\bfX_t|^p-|\bfX_t^\ast|^p\big|^{q/p}\big]\,,
\eeao
and the \rhs\ converges to zero in view of %\eqref{eq:couplingq}.
  \eqref{eq:couplingq2}.   
\par
By  a symmetry argument it suffices to bound the first term in $I_3$; the second one can be treated analogously. %The arguments are similar to those in Proposition \ref{pr:anticlpb}. 
Since  the anti-clustering condition is satisfied 
by $(|\bfX_t|)$  there exists $\theta>0$ such that
for every $x>0$,  we apply Theorem 3 of Segers~\cite{segers:2003} to obtain
\beao
\theta_n=\dfrac{\P(M_{r_n}^{|\bfX|}>x\,a_n)}{r_n\,\P(|\bfX|>x\,a_n)}\to\theta\,,\qquad \nto\,.
\eeao
Thus, for every $\vep>0$,
\beao
\lefteqn{k_n\P\Big(a_n^{-1}M_{r_n}^{|\bfX|}> x\,,a_n^{-1}M_{r_n}^{|\bfX^\ast|}\le x\,,\max_{1\le t\le r_n}|\bfX_t-\bfX_t^\ast|\le x\,\vep\, a_n\Big)}\\
&\le& k_n\,\big(\P\big(M_{r_n}^{|\bfX|}> x\, a_n\big)-
\P\big(M_{r_n}^{|\bfX|}> x\,(1+\vep)\, a_n\big)\big)\\
&=& \theta\,n\,\big(\P(|\bfX|> x\, a_n)-\P(|\bfX|> x(1+\vep)\, a_n)\big)+o(1)\\
&\to &\theta\,x^{-\a}\,(1-(1+\vep)^{-\a})\,,\qquad \nto\,,
\eeao
and the \rhs\ converges to $0$ as $\vep\downarrow 0$. 
Moreover, for every $x,\vep>0$,
\beao
\lim_{\nto}k_n\,\P\big(\max_{1\le t\le r_n}|\bfX_t-\bfX_t^\ast|> x\,\vep\, a_n\big)=0\,.
\eeao
This follows from \eqref{eq:couplingq} and an application of Markov's inequality of  the order 
$q$. Thus we proved that $\lim_{\nto} k_n\,I_3=0$. The proof is finished.
\end{proof}

\subsection{Iterated random Lipschitz functions}
Assume $(\bfX_t)$ is the solution of a system of iterated random functions: there exists a function $g$ and an iid sequence $(\vep_t)$ \st
\beam\label{eq:irlf}
\bfX_{t}=g(\bfX_{t-1},\vep_{t}),\qquad t\in \Z\,.
\eeam
We assume that this system is {\em contractive}, i.e., there exist
 $0<\rho<1$ and $q>0$ such that
\beam\label{eq:contr}
\E[|g(\bfx_0,\vep_0)- g(\bfx_1,\vep_0)|^q]\le \rho\,|\bfx_0-\bfx_1|^q\,,\qquad \bfx_0,\bfx_1\in \R^d\,. 
\eeam
If there exists $\bfx\in\R^d$ such that 
\beam\label{eq:condmom}
\E[|g(\bfx,\vep_0)|^q]<\infty\,,
\eeam 
the fixed point theorem in the complete space $L^q$ ensures 
the existence of a unique stationary solution $(\bfX_t)$ 
which admits finite moments of order $q>0$.
\par
A coupled version $(\bfX_t^\ast)$ is easily obtained as follows:
\beao
\bfX_t^\ast = \begin{cases}g(\bfX_{t-1}^\ast,\epsilon_t),&t\ge1\,,\\
g(\bfX_{t-1}^\ast,\epsilon_t'),& t\le 0\,,\end{cases}
\eeao
where $(\epsilon_t')$ is an independent copy of $(\epsilon_t)$. We have
\bpr\label{prop:IFS}
Assume that the unique stationary solution $(\bfX_t)$ of the
iterated random function system \eqref{eq:irlf} is \regvary\ 
with index $\a>0$ and satisfies \eqref{eq:contr}, 
\eqref{eq:condmom} for some $\rho\in (0,1)$ and $q<\a\wedge 1$. Then the conditions of Proposition \ref{prop:coupledmixed}
are satisfied.
\epr
\begin{proof} 
A recursive argument yields for some constant $c>0$,
\beao
\E\big[|\bfX_t-\bfX_t^\ast|^q]\le c \, \rho^t \,,\qquad t\ge 1\,.
\eeao
Then (2) in Proposition~\ref{prop:coupledmixed}  follows  choosing $C>0$ sufficiently large in the intermediate \seq\ $\ell_n=[C\log n]$ such that
\beao
 k_n\,a_n^{-q}\sum_{t=\ell_n}^{r_n}\big(\E\big[|\bfX_t-\bfX_t^\ast|^q\big]\big)^{(1/p)\vee 1} \le c \rho^{\ell_n}\dfrac{n}{r_na_n^q}\le c \dfrac{n^{C \log \rho+1-\delta}}{a_n^q}\to 0\,, \qquad \nto\,.
\eeao
It remains to verify (3) in Proposition~\ref{prop:coupledmixed} .We also have
\beao
\E\big[\big(\big|a_n^{-1} \big(\bfX_t - \bfX_t^*\big)\big|^{q}\wedge 1\big) \mid \bfX_{t-1},\bfX_{t-1}^\ast \big] &\le& \E\big[| a_n^{-1}(\bfX_t - \bfX_t^*)|^{q} \mid \bfX_{t-1},\bfX_{t-1}^\ast\big] \wedge 1\\
&\le&\big(\rho\,  \big| a_n^{-1}\big(\bfX_{t-1} - \bfX_{t-1}^*\big)\big|^{q}) \wedge 1 \,.
\eeao
Thus, using the filtration $\mathcal F_t = \sigma(\epsilon_t,\ldots,\epsilon_1,\bfX_0,\bfX_0')$, $t\ge 1$, and the Markov property we obtain
\beao
\lefteqn{\E\big[\big(\big|a_n^{-1} \big(\bfX_t - \bfX_t^*\big)\big|^{q}\wedge 1\big) \big(\big|a_n^{-1}\bfX_0\big|^{q }\wedge 1\big)\big]}\\
&\le&
\E\big[\E\big[\big(\big|a_n^{-1} \big(\bfX_t - \bfX_t^*\big)\big|^{q}\wedge 1\big)\mid \mathcal F_{t-1} \big]\big(\big|a_n^{-1}\bfX_0\big|^{q }\wedge 1\big)\big]\\
&\le&\E\big[\big(\big(\rho\,  \big| a_n^{-1}\big(\bfX_{t-1} - \bfX_{t-1}^*\big)\big|^{q}) \wedge 1\big)\big(\big|a_n^{-1}\bfX_0\big|^{q }\wedge 1\big)\big]\\
&\le&\E\big[\big(\big(\rho^t\,  \big| a_n^{-1}\big(\bfX_{0} - \bfX_{0}^*\big)\big|^{q}) \wedge 1\big)\big(\big|a_n^{-1}\bfX_0\big|^{q }\wedge 1\big)\big]\,,
\eeao
again using a recursive argument in the last step. Observing that
$\bfX_0$ and $\bfX_0^\ast$ are independent, we obtain
\beao
\lefteqn{\E\big[\big(\big|a_n^{-1} \big(\bfX_t - 
\bfX_t^*\big)\big|^{q}\wedge 1\big) \big(\big|a_n^{-1}\bfX_0\big|^{q }\wedge 1\big)\big]}\\&\le& \E\big[\big(\,  |\rho^{t/q} a_n^{-1}\bfX_{0}|^{q}\wedge 1\big) \big(\big|a_n^{-1}\bfX_0\big|^{q }\wedge 1\big)\big]+\rho^t\,a_n^{-2q}\E\big[\big|\bfX_0\big|^{q}\big]^2\\
&\le&\Big(\E\big[\big(\,  |\rho^{t/q} a_n^{-1}\bfX_{0}|^{q}\wedge 1\big)^2\big]\,\E\big[ \big(\big|a_n^{-1}\bfX_0\big|^{q }\wedge 1\big)^2\big]\Big)^{1/2}+\rho^t\,a_n^{-2q}\big(\E\big[\big|\bfX_0\big|^{q}\big]\big)^2\,.
\eeao
By Karamata's theorem there exists a positive constant $c>0$ such that 
\beao
\E\big[\big(\,  |\rho^{t/q} a_n^{-1}\bfX_{0}|^{q}\wedge 1\big)^2\big]&\le& c \rho^{\a t/q}\P(|X_0|>a_n)\\
\E\big[\big(\,  |a_n^{-1}\bfX_{0}|^{q}\wedge 1\big)^2\big]&\le& c\,\P(|X_0|>a_n)\,.
\eeao
For $q>\a/2$ we achieve
\beao
n\E\big[\big(\big|a_n^{-1} \big(\bfX_t - \bfX_t^*\big)\big|^{q}\wedge 1\big) \big(\big|a_n^{-1}\bfX_0\big|^{q }\wedge 1\big)\big]&\le& c\rho^{\a t/(2q)}+o(1)\,.
\eeao
Then (3) in  Proposition~\ref{prop:coupledmixed} and the desired result follow.
\end{proof}
\subsection{Examples}
In this section we consider two examples of \regvary\ stationary \ts : 
an autoregressive process or order 1 (AR(1)) and the solution to
an affine \sre\ (SRE). 
\bexam\rm {\bf A \regvary\ AR(1) process.}
We consider the causal stationary solution of the AR(1) equations 
$X_t=\varphi X_{t-1}+Z_t$, $t\in\bbz$, for some $\varphi\in (-1,1)\backslash\{0\}$ 
and an iid \regvary\ noise \seq\ $(Z_t)$. This means that a generic element 
$Z$ satisfies the tail balance condition, for $q_\pm\ge 0$ \st\ $q_++q_-=1$,
\beao
\dfrac{\P(\pm Z>x)}{\P(|Z|>x)}\to q_{\pm}\,,\qquad \xto\,.
\eeao
Then $(Z_t)$  has spectral tail process $\P(\Theta_0^Z=\pm 1)=q_{\pm}$, $\Theta_t^Z=0$, $t\ne 0$.
It is well known (e.g. Kulik and Soulier \cite{kulik:soulier:2020})
that $(X_t)$ is \regvary\ with index $\a$ and spectral tail process 
\beao
\Theta_t&=&
\Theta_Z\, \sign(\varphi^{J+t})|\varphi|^t\,\1(J+t\ge 0)=\Theta_Z\, \sign(\varphi^{J})\varphi^t\,\1(t\ge -J)=\Theta_0\,\varphi^t\,\1(t\ge -J)\,,\qquad
t\in\bbz\,,
\eeao
where $\P (\Theta_0=\pm 1)=p_\pm$, $J$ and $(\Theta_Z)$ are independent, and
\beao
p_\pm&=&q_\pm \,\1(\varphi>0)+\dfrac{q_\pm+q_\mp|\varphi|^\a}{1+|\varphi|^\a}\,\1(\varphi<0)\,\\
\P(J=j)& =& |\varphi|^{\a\,j}\,(1-|\varphi|^\a)\,,\qquad j=0,1,\ldots\,.
\eeao
The {\em forward spectral tail process} is given by
$\Theta_t= \Theta_0\,\varphi^t$, $t\ge 0$, and the spectral cluster process $(Q_t)$
by $Q_t=\Theta_t/\|\Theta\|_\a= \Theta_t\,(1-|\varphi|^\a)^{1/\a}$, $t\in \bbz$, and 
the extremal index by $\theta_{|X|}=1-|\varphi|^\a$.
\par
We observe that $X_t=g(X_{t-1},Z_t)=\varphi\, X_{t-1}+Z_t$, $t\in\bbz$, 
constitute a 
contractive iterative \fct\ system: for $q<\a$,
\beao
\E\big[|g(x_0,Z_0)-g(x_1,Z_0)|^q\big]= |\varphi|^q\,|x_0-x_1|^q\,.
\eeao
Hence the conditions of Proposition~\ref{prop:IFS} are satisfied and all 
 mixing and anti-clustering conditions needed for the results in this paper 
hold.
\eexam
\bexam \rm {\bf A \regvary\ solution to a SRE.} 
We consider the causal solution to an affine SRE $X_t=A_tX_{t-1}+B_t$, $t\in\bbz$, where $(A_t,B_t)$, $t\in\bbz$, is an iid $\bbr^2$-valued \seq .
We further assume that a generic element $(A,B)$ of this \seq\ satisfies
the following conditions: (i) there exists $\a>0$ \st\
$\E[|A|^\a]=1$, $\E[|A|^\a\log^+ |A|]<\infty$, $\E[|B|^\a]<\infty$, (ii)  the conditional law of 
$\log |A|$ given $\{A\ne 0\}$ is non-arithmetic, (iii) $\P(A\,x+B=x)<1$ for every 
$x\in\bbr$. Then there exists an a.s. unique causal solution to the SRE
with the property $\P(\pm X_0>x)\sim c_\pm \,x^{-\a}$ as $\xto$ 
for constants $c_\pm$ \st\ $c_++c_->0$. This follows from classical
Kesten-Goldie theory; cf. Theorem 2.4.7 in Buraczewski et al. \cite{buraczewski:damek:mikosch:2016}. The \seq\ $(X_t)$ is \regvary\ with index $\a$,
forward spectral tail process $\Theta_t=\Theta_0\,A_1\cdots A_t$, $t\ge 0$,
where $\P(\Theta=\pm 1)=c_\pm/(c_++c_-)$ and extremal index 
$\theta_{|X|}= \P\big[\big(1-\sup_{t\ge 1}|A_1\cdots A_t|^\a\big)_+\big]$;
see Basrak and Segers \cite{basrak:segers:2009}.
\par
We observe that $X_t=g(X_{t-1},(A_t,B_t))=A_tX_{t-1}+B_t$, $t\in\bbz$,
constitute a contractive iterative function system: for $0<q<\a$, with $\rho=\E[|A_0|^q]$,
$\E[ |g(x_0,(A_0,B_0))-g(x_1,(A_0,B_0))|^q]= \rho\,|x_0-x_1|^q\,.$
The value $\rho<1$ by convexity of the \fct\ $f(q)=\E[|A_0|^q]$ and since
$f(\a)=1$. Hence the conditions of Proposition~\ref{prop:IFS} are satisfied and all 
 mixing and anti-clustering conditions needed for the results in this paper 
hold.
\eexam

\appendix
\section{Some auxiliary results}\setcounter{equation}{0}
\ble\label{lem:xmas29a}  
Consider  an $\bbr^d$-valued stationary \regvary\ \seq\ $(\bfX_t)$
with index $\a\in (1,2)$. If the anti-clustering condition \eqref{cond:acac} 
is satisfied then 
\beam\label{eq:xmas29e}
J:=\E\Big[\Big(\sum_{j= 0}^\infty|\bfTh_j|\Big)^{\a-1}\Big]<\infty\,.
\eeam
\ele
\begin{proof}    
By sub-additivity we have
\beao
J
&\le& \E\Big[\Big(\sum_{j=0}^\infty|\bfTh_j|\,\1(|\bfTh_j|\le1)\Big)^{\a-1}\Big]+\E\Big[\sum_{j=0}^\infty|\bfTh_j|^{\a-1}\,\1(|\bfTh_j|>1)\Big]
=:I_1+I_2\,.
\eeao
By  Jensen's inequality,  
\beao
I_1\le \Big(\E\Big[ \sum_{j=0}^\infty|\bfTh_j|\,\1(|\bfTh_j|\le1) \Big]\Big)^{\a-1} \,.
\eeao
We prove that the \rhs\ is finite by showing
\beam\label{eq:sumtrunc}
\sum_{j=0}^\infty\E\big[|\bfTh_{j}|\wedge  1 \big]<\infty\,.
\eeam
Using the decomposition
\beao
 |a_n^{-1}\bfX_j |\wedge  1 &= &   |a_n^{-1}\bfX_j|\,\1(|\bfX_j|\le a_n) 
+ \1(|\bfX_j|>a_n)\,,
\eeao
condition \eqref{cond:acac} implies
\beao
\lim_{k\to\infty} \limsup_{\nto}\sum_{j=k}^{r_n}  \,n\, \,\E\big[(|a_n^{-1}\bfX_{j} |\wedge  1)\,\1(|  \bfX_{0}|>
a_n)\big]=0\,.
\eeao
 Hence for every 
$\vep>0$ there exists an integer $k_0$ sufficiently large \st 
\beao
\limsup_{\nto}\sum_{j=k}^{k+h}  \,n\, \,\E\big[(|a_n^{-1}\bfX_{j} |\wedge  1)\,\1(|  \bfX_{0}|>
a_n)\big]\le \vep\,,\qquad k\ge k_0\,, \qquad h\ge 0\,.
\eeao
Using the  \regvar\ of $(\bfX_t)$ with tail process $(Y\,\bfTh_t)$, 
we can determine the limit of each summand
\beao
 \lim_{\nto}  n\, \,\E\big[(|a_n^{-1}\bfX_{j} |\wedge  1)\,\1(|  \bfX_{0}|>
a_n)\big]&=& \lim_{\nto}   \E\big[(|a_n^{-1}\bfX_{j} |\wedge  1)\; \big| \; |\bfX_{0}|>a_n\big]\\
&=&    \E\big[|Y\,\bfTh_{j}|\wedge  1 \big]\,.
\eeao
Thus the Cauchy criterion for the \seq\ 
$(\sum_{j=0}^k\E\big[|Y\bfTh_{j}  |\wedge  1 \big])_{k\ge 0}$ holds. 
Since $Y>1$ a.s. \eqref{eq:sumtrunc} follows.
\par
It remains to show that  $I_2<\infty$.
By stationarity we can show similarly to \eqref{eq:sumtrunc} that
\beao
 \sum_{j=0}^{\infty}
    \E[|\bfTh_{-j}|\land 1] < \infty. 
\eeao 
By the time-change formula  \eqref{eq:tcp} on p.~\pageref{eq:tcp} we obtain that 
\beao
\E[|\bfTh_{-j}|\land 1]&=&\E[|\bfTh_{-j}|\land 1\mid \bfTh_{-j}\neq {\bf0}]\P(\bfTh_{-j}\neq {\bf0})\\
&=&\E[|\bfTh_{-j}|\land |\bfTh_0|\mid \bfTh_{-j}\neq {\bf0}]\E\big[|\bfTh_{j}|^\a\big]\\
&=& \E[|\bfTh_j|^\alpha\,( |\bfTh_{j}|^{-1}\land 1)]\,.
\eeao
We conclude that
\beao
    \infty &>&\sum_{j=0}^{\infty} \E[|\bfTh_{-j}|\land 1]= \sum_{j=0}^{\infty} \E[|\bfTh_j|^\alpha\,( |\bfTh_{j}|^{-1}\land 1)] \\
    &=& \sum_{j=0}^{\infty} \E[ |\bfTh_{j}|^{\alpha-1}\land |\bfTh_j|^\alpha]>
    \sum_{j=0}^{\infty} \E[ |\bfTh_{j}|^{\alpha-1}\1(|\bfTh_j| > 1)],
\eeao
implying $I_2<\infty$.

\end{proof}
\section{Proof of Lemma \ref{lem:domratio}}\label{app:domratio}\setcounter{equation}{0}
\begin{proof}
For $1<\a<2$, $\bfxi_\a$ is  $\a$-stable and $\E[|\bfxi_\a|]<\infty$ which implies the statement of Lemma  \ref{lem:domratio}.

Next we consider the case $0<\a<1$. We use the series representation \eqref{eq:april10a}:
\beao
|\bfxi_\a|=\Big|\sum_{i=1}^\infty \Gamma_i^{-1/\a} \sum_{j\in\Z} \bfQ_{ij}\Big|\le \sum_{i=1}^\infty \Gamma_i^{-1/\a} \sum_{j\in\Z} |\bfQ_{ij}|=:\xi_\a'\,.
\eeao
The \rhs\ represents a positive $\a$-stable random variable $\xi_\a'$.
Indeed,  $(\sum_{j\in\Z} |\bfQ_{ij}|)_{i\in\Z}$ is a \seq\ of iid \rv s satisfying $\E[(\sum_{j\in\Z} |\bfQ_{j}|)^\a]<\infty$ since $\|\bfQ\|_\a=1$ by definition. 
We apply Fubini's theorem for positive \rv s and obtain 
\beao
\E\Big[\dfrac{\xi_{\a}'}{(\zeta_{\a,p}^p)^{1/p}}\Big]&=&\dfrac  p{\Gamma(1/p)}\int_0^\infty \E\big[\xi_{\a}'\ex^{-\lambda^p\zeta_{\a,p}^p}\big]d\lambda\,.
\eeao
We will show that the integrand coincides with 
\beao
\lim_{x\to 0^+}-\dfrac{\partial \Psi_{\xi_\a',\zeta_{\a,p}^p}(x,\lambda^p)}{\partial x}\,,\qquad  \Psi_{\xi_\a',\zeta_{\a,p}^p}(x,\lambda^p):=\E\big[\ex^{-x\, \xi_\a'-\lambda^p\, \zeta_{\a,p}^p}\big]\,.
\eeao
From the expression of the \chf\ - Laplace transform of Theorem~\ref{pr:prophybridch} we have
\beao
\Psi_{\xi_\a',\zeta_{\a,p}^p}(x,\lambda^p)&=&\exp\Big(\dint_0^\infty\E\Big[\ex^{-y\,x\, \| \bfQ\|_1-y^p\la \|\bfQ\|_p^p}-1\Big]d(-y^{-\a})\Big)\\
-\dfrac{\partial \Psi_{\xi_\a',\zeta_{\a,p}^p}(x,\lambda^p)}{\partial x}&=&\int_0^\infty \E\Big[y\,\|\bfQ\|_1
\ex^{-xy\|\bfQ\|_1-(\la y)^p\|\bfQ\|_p^p} \Big]
d(-y^{-\a})\;\E\big[\ex^{-x\xi_\a'-\lambda^p\zeta_{\a,p}^p}\big]\\
&=&\E\Big[\|\bfQ\|_p^\a\dfrac{\|\bfQ\|_1}{\|\bfQ\|_p}\int_0^\infty y\,
\ex^{-xy\|\bfQ\|_1/\|\bfQ\|_p-(\la y)^p} 
d(-y^{-\a})\Big]\;\E\big[\ex^{-x\xi_\a'-\lambda^p\zeta_{\a,p}^p}\big]\\
&\to&\dfrac\a p \Gamma((1-\a)/p)\la^{\a-1}\E\Big[\|\bfQ\|_p^\a\dfrac{\|\bfQ\|_1}{\|\bfQ\|_p}\Big]\;\E\big[\ex^{-\lambda^p\zeta_{\a,p}^p}\big]\,,\qquad x\to 0^+\,,
\eeao
where we exploit \eqref{eq:gammaformula} in the last step. By monotone convergence we also have 
\beao
\lim_{x\to 0^+}-\dfrac{\partial \Psi_{\xi_\a',\zeta_{\a,p}^p}(x,\lambda^p)}{\partial x}=\lim_{x\to 0^+}\E\big[ \xi_\a'\ex^{-x\, \xi_\a'-\lambda^p\, \zeta_{\a,p}^p}\big]=\E\big[ \xi_\a'\ex^{-\lambda^p\, \zeta_{\a,p}^p}\big]
\eeao
since the limit exists. Thus $\E\big[ \xi_\a'\ex^{-\lambda^p\, \zeta_{\a,p}^p}\big]<\infty$ and we conclude the proof of Lemma \ref{lem:domratio} using the domination $|\bfxi_\a|\le \xi_\a'$.
\end{proof}
\section{Proof of the integrability of $\bfR_\a$}\label{app:domratiobis}\setcounter{equation}{0}
We proceed as in Appendix \ref{app:domratio} dominating
$|\bfR_\a|\le\xi_\a'/\eta_\a$ for $\a\in(0,1)$. The integrability of $\bfR_\a$, $\a\in(1,2)$, follows easily from $\E[|\bfxi_\a|]<\infty$. We introduce the hybrid Laplace transform
\beao
 \Psi_{\xi_\a',\eta_\a}(u,x):=\E\big[\ex^{-u\, \xi_\a'}\1(\eta_\a\le x)\big]\,,\qquad u,x>0.
\eeao
From the expression of the hybrid \chf\  of Theorem~\ref{th:joint}  we have
\beao
 \Psi_{\xi_\a',\eta_\a}(u,x)= \exp\Big(\int_0^\infty\E\big[\ex^{-yu\|\bfQ\|_1}\1(y\|\bfQ\|_\infty\le x)-1\big]d(-y^{-\a})\Big)\,,\qquad u,x>0\,.
 \eeao
Then by monotone convergence we obtain 
\beao
\lim_{u\to 0^+}-\dfrac{\partial \Psi_{\xi_\a',\eta_{\a}}(u,x)}{\partial u}&=&\lim_{u\to 0^+}\int_0^\infty\E\big[y\|\bfQ\|_1\ex^{-yu\|\bfQ\|_1}\1(y\|\bfQ\|_\infty\le x)]d(-y^{-\a}) \Psi_{\xi_\a',\eta_\a}(u,x)\\
&=&\a\E\Big[\|\bfQ\|_1\int_0^{x/\|\bfQ\|_\infty}y^{-\a}dy\Big] \Psi_{\xi_\a',\eta_\a}(0,x)\\
&=&\dfrac\a{1-\a}\E\Big[\dfrac{\|\bfQ\|_1}{\|\bfQ\|_\infty^{1-\a}}\Big]x^{1-\a}\exp(-\theta_{|\bfX|}x^{-\a})\,.
\eeao
The limit is finite since we show the integrability of $\sum_{t\in\Z}|\wt \bfQ_t|$ in Remark \ref{rem:existence:sumQ} and
\beao
\E\Big[\dfrac{\|\bfQ\|_1}{\|\bfQ\|_\infty^{1-\a}}\Big]=\E\Big[\|\bfQ\|_\infty^\a\dfrac{\|\bfQ\|_1}{\|\bfQ\|_\infty}\Big]=\theta_{|\bfX|}\E[\|\wt \bfQ\|_1]<\infty\,.
\eeao
Then we proved the finiteness of 
\beao
\E[\xi_\a'\1(\eta_\a\le x)]=\lim_{u\to 0^+}-\dfrac{\partial \Psi_{\xi_\a',\eta_{\a}}(u,x)}{\partial u}=\dfrac\a{1-\a} \theta_{|\bfX|}\E[\|\wt \bfQ\|_1]x^{1-\a}\exp(-\theta_{|\bfX|}x^{-\a})\,,\qquad x>0\,.
\eeao
Applying Fubini's theorem we achieve
\beao
\E\Big[\dfrac{\xi_\a'}{\eta_\a}\Big]&=&\E\Big[\xi_\a'\int_0^\infty\1(\eta_\a\le x)\dfrac{dx}{x^2}\Big]\\
&=&\int_0^\infty\E[\xi_\a'\1(\eta_\a\le x)]\dfrac{dx}{x^2}\\
&=&\dfrac{\E[\|\wt \bfQ\|_1]}{1-\a}\int_0^\infty\a \theta_{|\bfX|}x^{1-\a}\exp(-\theta_{|\bfX|}x^{-\a})\dfrac{dx}{x^2}\\
&=&\dfrac{\E[\|\wt \bfQ\|_1]}{1-\a}\,.
\eeao
The \lhs\ term is finite which proves that $\bfR_\a$ is integrable for $0<\a<1$.

\end{document}